\documentclass[12pt,a4paper]{article}

\usepackage{amsmath}
\usepackage{amsfonts}
\usepackage{amssymb}

\numberwithin{equation}{section}

\newtheorem{thm}{Theorem}
\newtheorem{lem}{Lemma}
\newtheorem{prop}{Proposition}
\newtheorem{cor}{Corollary}
\newtheorem{rem}{Remark}

\newcommand*{\C}{\mathbb{C}}
\newcommand*{\R}{\mathbb{R}}
\newcommand*{\Q}{\mathbb{Q}}
\newcommand*{\Z}{\mathbb{Z}}
\newcommand*{\F}{\mathbb{F}_q}
\newcommand*{\A}{\mathbb{A}}

\begin{document}

\title{\bf The Riemann hypothesis for Weng's zeta function of ${\rm Sp}(4)$ over $\Q$
}
\author{Masatoshi Suzuki~\thanks{This work is supported by JSPS.}}

\maketitle

\begin{abstract}
As a generalization of the Dedekind zeta function,
Weng defined the high rank zeta functions and proved that
they have standard properties of zeta functions, namely,
meromorphic continuation, functional equation, and having only two simple poles.
The rank one zeta function is the Dedekind zeta function.
For the rank two case, the Riemann hypothesis is proved for a general number field.
Recently, he defined more general new zeta function associated to
a pair of reductive group and its maximal parabolic subgroup.
As well as high rank zeta functions,
the new zeta function satisfies standard properties of zeta functions.
In this paper, we prove that the Riemann hypothesis of Weng's zeta function
attached to the sympletic group of degree four.
This paper includes an appendix written by L. Weng,
in which he explains a general construction for zeta functions
associated to ${\rm Sp}(2n)$.
\end{abstract}

\section{Introduction}

Let $C$ be a nonsingular projective algebraic curve over a finite field $\F$.
The zeta function $\zeta_C(s)$ is defined by the Euler product
$\zeta_C(s)=\prod_{P \in C_0}(1-q^{-s \deg(P)})^{-1}$,
where $C_0$ is the set of all closed points of $C$.
The theory of $\zeta_C(s)$ is one of the most beautiful and successful one in Number theory.
A lot of nice theory was established by modeling on the theory of $\zeta_C(s)$.
Recently, Lin Weng defined a class of general zeta functions starting from the following formula of $\zeta_C(s)$:
\begin{equation} \label{101}
\zeta_C(s) = \sum_{[L] \in {\rm Pic}(C)} \frac{q^{h^0(L)}-1}{q-1} q^{-s \deg(L)},
\quad \Re(s)>1,
\end{equation}
where ${\rm Pic} \, (C)$ is the moduli space of isomorphic classes $[L]$ of line bundle $L$ on $C$,
and $h^0(L)$ is the dimension of the $0$-th cohomology $H^0(C,L)$ over $\F$.
First he generalized $\zeta_C(s)$ to the rank $r \geq 1$ zeta function $\zeta_{C,r}(s)$
by replacing ${\rm Pic} \, (C)$
by the moduli space of isomorphic classes of semi-stable vector bundles on $C$ of rank $r$.
The semi-stability (introduced by Mumford) is needed to get a reasonable structure of moduli space
which is reflected to the standard properties of zeta function $\zeta_{C,r}(s)$.

On the other hand, for the number field case,
Iwasawa's interpretation of the Dedekind zeta function
allows us to get an arithmetic analogue of \eqref{101}
for the Dedekind zeta function.
Let $F$ be a number fields with discriminant $\Delta_F$,
${\A}_F^\times$ be the idele group of $F$
and $w_F$ be the number of roots of unity in $F$.
We denote by $\xi_F(s)$ the completed Dedekind zeta function of $F$.
Any idele $a = (a_v)_v \in {\A}_F^\times$
defines an Arakelov line bundle $L_a$ on $\overline{X}_F={\rm Spec} \, {\cal O}_F \cup X_{F,\infty}$,
and any Arakelov line bundle is obtained by this manner.
(More precisely, the finite part of $(a_v)_v$ defines a invertible sheaf on ${\rm Spec}\,{\cal O}_F$,
and the infinite part of $(a_v)_v$ defines a Hermitian metric on $X_{F,\infty}$.)
Moreover we have the topological group isomorphism ${\A}_F^\times/F^\ast \simeq {\rm Pic}(F)$,
where ${\rm Pic}(F)$ is the moduli space of isomorphic classes of Arakelov line bundles on $\overline{X}_F$.
Throughout this relation between ideles and Arakelov line bundles,
Iwasawa's interpretation gives
\begin{equation} \label{102}
\xi_F(s) = w_F^{-1} |\Delta_F|^{\frac{s}{2}} \int_{{\rm Pic}(F)} (e^{h^0(L)}-1) e^{-s \deg(L)} d\mu(L),
\quad \Re(s)>1,
\end{equation}
where $\deg(L)$ is the arithmetic degree of an Arakelov line bundle $L$
and $d\mu$ is a certain natural measure on ${\rm Pic}(F)$.
In \eqref{102}, the definition of the arithmetic (dimension of) cohomology $h^0(L)=h^0( \overline{X_F},L)$
is naturally given on the line of Iwasawa's interpretation (see \cite[B.2.1]{We1},
and also ~\cite{GS}),
and Tate's Fourier analysis on ${\A}_F$ gives the Riemann-Roch theorem
as a consequence of the Pontragin duality and the Poisson summation formula.
The Pontragin duality is interpreted as the Serre duality under suitable definition of $h^1$,
see \cite[B.2.3.3, B.2.3.4]{We1}, and also \cite[section 1.8]{We2}.

Standing on the arithmetic-geometric formula \eqref{102}
and Tate's Fourier analysis on ${\A}_F$,
Weng defined the (completed) rank $r$ zeta function $\xi_{F,r}(s)$ of $F$ by
\begin{equation} \label{103}
\xi_{F,r}(s) = |\Delta_F|^{\frac{rs}{2}} \int_{{\cal M}_{F,r}} (e^{h^0(E)}-1) e^{-s \deg(E)} d\mu(E),
\quad \Re(s)>1,
\end{equation}
where ${\cal M}_{F,r}$ is the moduli space of isomorphic classes $[E]$ of
rank $r$ semi-stable Arakelov vector bundles $E$ on $\overline{X}_F$,
and $d\mu$ is its associated Tamagawa measure, and $\deg (E)$ is the arithmetic degree of $E$.
The arithmetic analogue of $h^0$ and $h^1$ are defined as a natural extension of the line bundle case
having in mind Tate's method (\cite[B.2.3.4]{We1}, see also \cite{Ar1}).
We have $\xi_{F,1}(s) = w_F \xi_{F}(s)$, since any Arakelov line bundles are semi-stable.
Thanks of the semi-stability the right-hand side of \eqref{103} converges absolutely for $\Re(s)>1$.
The arithmetic-geometric Riemann-Roch theorem is proved
by using the Fourier analysis on ${\A}_F^n$ following Tate (\cite[B.2.3.4]{We1}).
As a consequence of the Serre duality and the Riemann-Roch theorem,
the rank $r$ zeta function $\xi_{F,r}(s)$ is continued meromorphically to $\C$
with only two simple poles at $s=0,1$,
and satisfies the standard functional equation $\xi_{F,r}(s)=\xi_{F,r}(1-s)$ (\cite[B.2.4.2]{We1}).
Hence the Riemann hypothesis for $\xi_{F,r}(s)$ is stated
that all zeros of $\xi_{F,r}(s)$ lie on the central line $\Re(s)=1/2$.

Remarkable fact for $\xi_{F,r}(s)$ is that the Riemann hypothesis for the rank $2$ zeta function $\xi_{F,2}(s)$
is \emph{proved} for any number field $F$ (\cite{LaSu},~\cite{We2}, see also \cite{Ha}).
Moreover all zeros of $\xi_{{\Q},2}(s)$ are simple (\cite{LaSu},~\cite{Ki}).
Now it is expected that the Riemann hypothesis holds for any $r \geq 1$ (\cite{We3}).
Note that the case $r=1$ is the Riemann hypothesis of the Dedekind zeta function.
The study of high rank zeta function $\xi_{F,r}(s)$ lead us to more general zeta functions.
%

In general,  $\xi_{{\Q},r}(s)$ is expressed as an integral of the Eisenstein series (\cite{We3}).
For example, the rank two case is
\begin{equation} \label{104}
\xi_{{\Q},2}(s) = \int_{{\cal M}_{{\Q},2}[1]} \widehat{E}(z,s) d\mu(z),
\end{equation}
where ${\cal M}_{{\Q},2}[1] \simeq D_0$, $D_T=\{ z=x+iy \,|\, |x| \leq 1/2, ~ 0< y \leq \exp(T), ~|z| \geq 1   \} $ and
\[
\widehat{E}(z,s) = \xi(2s) \!\!\!   \sum_{\gamma \in P(\Z) \backslash {\rm SL}(2,\Z)} \!\!\!  \Im(\gamma z)^s,
\quad (\Re(s)>1), \quad  P = B = \left\{  \begin{pmatrix} 1 & \ast \\ 0 & 1 \end{pmatrix} \right\}.
\]
Recall that the group theoretic description of the upper half plane
\[
D_0 \hookrightarrow {\cal H} \simeq {\rm SL}(2,{\Z}) \backslash {\rm SL}(2,{\R}) \slash {\rm SO}(2)
\]
and the fact that $P$ is the maximal parabolic subgroup of ${\rm SL}(2)$.
These lead us to more general new zeta functions attached to a pair $(G,P)$ of a reductive group $G$
and a maximal parabolic subgroup $P$ of $G$ motivated by \eqref{104}.

First important point for the generalization is the relation between the geometric truncation and the analytic truncation.
In the rank two case, define the analytic truncation $\wedge^T \widehat{E}(z,s)$ by
\[
\wedge^T \widehat{E}(z,s) =
\begin{cases}
\widehat{E}(z,s), & \text{if $y \leq \exp(T)$, $z \in D_\infty$},\\
\widehat{E}(z,s) - a_0(y,s), & \text{if $y > \exp(T)$, $z \in D_\infty$}.
\end{cases}
\]
where $a_0(y,s)$ is the constant term of the Fourier expansion of $\widehat{E}(z,s)$.
Then we have
\begin{equation} \label{Geo_Ana}
\int_{D_T} \widehat{E}(z,s) d\mu(z)
= \int_{{\rm SL}(2,{\Z}) \backslash {\rm SL}(2,{\R}) \slash {\rm SO}(2)} \wedge^T \widehat{E}(z,s) d\mu(z).
\end{equation}
Namely, the geometric truncation $D_T$ is equal to the analytic truncation $\wedge^T$ in the rank two case.
Taking $T=0$, we obtain the rank two zeta $\xi_{{\Q},2}(s)$.
This kind of equality holds in widely general situation if $T$ is sufficiently regular
(see A.2.2 of Appendix~\cite{We4}).

To explain the second point for the generalization,
we restrict as $G={\rm Sp}(4)$ and $F={\Q}$ for the simplicity.
See Appendix~\cite{We4} for the detailed theory of ${\rm Sp}(2n)$,
\cite{SW} for the zeta functions of the exceptional group $G_2$,
and \cite{We5} for the general theory of zeta functions associated to the pair
$(G,P)$ of reductive group and its maximal parabolic subgroup.

Let ${\frak S}_2$ be the Siegel upper half-space of degree $2$,
and let $\Gamma_2=\{\pm I_4 \} \backslash {\rm Sp}(4,\Z)$
be the Siegel modular group of degree $2$, where $I_4$ is the identity matrix of size $4$.
Any $Z \in {\frak S_2}$ is written as $Z=X+\sqrt{-1}Y$ with $X,Y \in {\rm Sp}(4,\R)$
such that $Y= \Im(Z)$ is positive definite.
For $Z \in {\frak S_2}$ and $\Re(s) \ll 0$,
the Siegel-Maa\ss~Eisenstein series is defined by
\begin{equation*}
E_2(Z,s)=\sum_{\bigl( \smallmatrix \ast & \ast \\ C & D \endsmallmatrix \bigr) \in {\frak P}_2 \backslash \Gamma_2}
\frac{ |Y|^{-s}}{\Vert CZ+D \Vert^{-2s}}
\end{equation*}
where ${\frak P}_2=\{ \bigl( \smallmatrix \ast & \ast \\ 0 & \ast \endsmallmatrix \bigr) \} \cap \Gamma_2 $,
$|Y|=\det \, Y$
and $\Vert CZ+D \Vert = |\det(CZ+D)|$.
Fix a minimal parabolic subgroup (a Borel subgroup) $P_0$ of ${\rm Sp}(4)$,
and a Levi component $M_{0}$ of $P_0$.
Let $X(M_{0})$ be the group of characters of $M_{0}$ defined over $\Q$,
$\Delta_0$ be a base for a root system,
${\frak a}_{0}={\rm Hom}_{\Z}(X(M_0),{\R})$
and ${\frak a}_{0}^+= \{ T \in {\frak a}_{0} \,|\, \langle \alpha,T\rangle > 0, ~\forall \alpha \in \Delta_0 \}$.
An element $T \in {\frak a}_{0}^+$ is called sufficiently regular
if $\langle \alpha, T \rangle>0$ is large enough for any $\alpha \in \Delta_0$.
For sufficiently regular $T \in {\frak a}_{0}^+$, we define
\begin{equation} \label{105}
Z_{{\rm Sp}(4),{\Q}}^T(s) := \int_{\Gamma_2 \backslash {\frak S}_2} \wedge^T E_2(Z;s)\,d\mu(Z),
\end{equation}
where $\wedge^T$ is Arthur's truncation operator with respect to $T$.
Let $\frak F(T)$ be the compact subset of $\Gamma_2 \backslash {\frak S}_2$
whose characteristic function is given by $\wedge^T \mathbf{1}$.
Then we have
\begin{equation} \label{period_0}
\aligned
Z_{{\rm Sp}(4),{\Q}}^T(s)
 = \int_{\Gamma_2 \backslash {\frak S}_2} \wedge^T E_2(Z;s)\,d\mu(Z)
 = \int_{\frak F(T)}  E_2(Z;s) \,d\mu(Z).
\endaligned
\end{equation}
for every sufficiently regular $T \in {\frak a}_{0}^+$ by Corollary 2.2.1 of Appendix \cite[A.2.2, A.2.3]{We4}.
Let recall the formula \eqref{104} of high rank zeta functions $\xi_{{\Q},r}(s)$.
Then we understand that $Z_{{\rm Sp}(4),{\Q}}^{T=0}(s)$ is a natural analogue of $\xi_{{\Q},r}(s)$,
if we can give a reasonable arithmetic or geometric meaning of the set $\frak F(0)$ (\cite[\S A.2.7]{We4}).
Here we have two problems:
\begin{enumerate}
\item[(1)] How to calculate the integral in \eqref{period_0}?
\item[(2)] Can we take $T=0$ in \eqref{period_0}?
\end{enumerate}
To attack these problems, we note the formula
\begin{equation*}
E_2(Z,s) = \underset{z=s+\frac{1}{2}}{\mathrm{Res}} \Big(E\Big({\bf 1} \,; z,s-\frac{1}{2}\,; Z \Big)\Big),
\end{equation*}
where $E(\varphi; z_1, z_2 \, ; Z)$ is the Langlands-Eisenstein series
associated to the Borel subgroup $P_0$, $(z_1,z_2) \in X(M_0) \otimes {\C} \simeq {\C}^2$,
$\varphi$ is cuspidal at the level of Borel, and $\mathbf{1}$ is the constant function one
(see A.2.5 of Appendix \cite{We4} for the definition of $E(\varphi; z_1, z_2 \, ; Z)$,
and \cite[section 2]{SW}], or better,
\cite[sections 2,3]{We3} for the general theory of this one).
Then we have
\begin{equation} \label{period_1}
Z_{{\rm Sp}(4),{\Q}}^T(s)
 = \int_{\frak F(T)} \underset{z=s+\frac{1}{2}}{\mathrm{Res}} \Big(E\Big({\bf 1} \,; z,s-\frac{1}{2}\,;Y\Big)\Big) \,d\mu(Z).
\end{equation}
Now we \emph{assume} that
\begin{enumerate}
\item[($\star$)] we can exchange the integration and the taking the residue in \eqref{period_1}.
\end{enumerate}
Then we have
\begin{equation} \label{period_2}
Z_{{\rm Sp}(4),{\Q}}^T(s)
 = \underset{z=s+\frac{1}{2}}{\mathrm{Res}} \int_{\frak F(T)}  \Big(E\Big({\bf 1} \,; z,s-\frac{1}{2}\,;Y\Big)\Big) \,d\mu(Z).
\end{equation}
The integral of the right-hand side is calculated \emph{explicitly}
in terms of the Weyl group and the Riemann zeta function by \cite[Corollary 17]{JaLaRo}.
Moreover, in that explicit description, we can take $T=0$! (see A.2.6 of Appendix ~\cite{We4}.)

Unfortunately, assumption ($\star$) is \emph{not} allowed in general,
even if $G={\rm SL}(n)~(n \geq 3)$.
In other words, we can not define a \emph{direct} analogue
$Z_{{\rm Sp}(4),{\Q}}^{T=0}(s)$ of $\xi_{{\Q},r}(s)$ at present
(more precisely, see A.2.7 of Appendix \cite{We4}).
Here we turn the consideration.
We define the \emph{new} zeta function by using the right-hand side of \eqref{period_2}.
This is the second point of the generalization mentioned above.
We define the new zeta function $\xi_{{\rm Sp}(4),{\Q}}(s)$ by
\begin{equation} \label{106}
\xi_{{\rm Sp}(4),{\Q}}(s) := \xi(2)\xi(s+1)\xi(2s) \cdot
\underset{z=s-\frac{1}{2}}{\mathrm{Res}}
\left[\,
\int_{\frak F(T)}  \Big(E\Big({\bf 1} \,; z,s-\frac{3}{2}\,;Y\Big)\Big) \,d\mu(Z)
\right]_{T=0}.
\end{equation}
Here the operation $T=0$ is justified via the explicit calculation of the integral
by using  \cite[Corollary 17]{JaLaRo} (see A.2.7 of Appendix~\cite{We4}),
the factor $\xi(2)\xi(s+1)\xi(2s)$ is to clearance the denominator of that explicit form,
and the change $s-\frac{1}{2} \mapsto s-\frac{3}{2}$ is to normalize the form of functional equation.
Zeta function of ${\rm Sp}(4)$ or more general reductive group $G$ is defined by a similar way.
As well as the high rank zeta function $\xi_{{\Q},r}(s)$,
the above new zeta functions satisfies a standard properties of zeta functions,
even if $G \not={\rm SL}(n), ~{\rm Sp}(2n)$.
We conjectured the functional equation for zeta functions of $(G,P)$,
but only verify it for a few concrete examples.
This is in fact the next most important yet doable
question in this direction.

As for ${\rm Sp}(4)$, $\xi_{{\rm Sp}(4),{\Q}}(s)$ is calculated explicitly as
\begin{equation} \label{107}
\begin{aligned}
\xi_{{\rm Sp}(4),{\Q}}(s)=
&~\frac{1}{s-2}\xi(2)\cdot\xi(s+1)\xi(2s)-\frac{1}{s+1}\xi(2)\cdot\xi(s-1)\xi(2s-1)\\
&-\frac{1}{2s-2}\cdot \xi(s+1)\xi(2s)+\frac{1}{2s}\cdot\xi(s-1) \xi(2s-1)\\
&-\frac{1}{(2s-2)(s+1)}\cdot\xi(s)\xi(2s)-\frac{1}{(2s)(s-2)}\cdot\xi(s)\xi(2s-1),
\end{aligned}
\end{equation}
where $\xi(s)=\pi^{-s/2}\Gamma(s/2)\zeta(s)$ is the completed Riemann zeta function
(see A.3.1 of Appendix~\cite{We4}).
By formula \eqref{107}, $\xi_{{\rm Sp}(4),{\Q}}(s)$ is continued meromorphically to $\C$,
satisfies the canonical functional equation,
\begin{equation}
\xi_{{\rm Sp}(4),{\Q}}(s) = \xi_{{\rm Sp}(4),{\Q}}(1-s)
\end{equation}
and  all poles of $\xi_{{\rm Sp}(4),{\Q}}(s)$ are four simple poles $s=0,1$ and $s=-1,2$.
For the zeros of $\xi_{{\rm Sp}(4),{\Q}}(s)$, we have the following significant result.
\begin{thm}[RH for ${\rm Sp}(4)/{\Q}$] \label{thm_01}
All zeros of $\xi_{{\rm Sp}(4),{\Q}}(s)$ lie on the line $Re(s)=1/2$.
\end{thm}
The entire function
\begin{equation} \label{108}
Z(s) := 4 s^2 (s-1)^2 \cdot (s+1)(2s-1)(s-2) \cdot \xi_{{\rm Sp}(4),{\Q}}(s).
\end{equation}
is more useful than $\xi_{{\rm Sp}(4),{\Q}}(s)$ itself for the proof of Theorem \ref{thm_01}.
We have
\begin{equation} \label{109}
\aligned
Z(s) & = (s-1)(As-A+1) \cdot \chi(s+1) \chi(2s)
   - (s-2) \cdot \chi(s)  \chi(2s)\\
& \quad - s (As-1) \cdot \chi(s-1) \chi(2s-1)
  - (s+1) \cdot \chi(s)  \chi(2s-1),
\endaligned
\end{equation}
where
\begin{equation}
A = 2\xi(2)-1=\pi/3-1>0, \quad
\chi(s) = s(s-1)\xi(s).
\end{equation}
Note that $Z(s)$ has real zeros at $s=0,1$ and $s=1/2$,
because $s=0,1$ are simple poles of $\xi_{{\rm Sp}(4),{\Q}}(s)$
and $s=1/2$ is its regular point.
Now Theorem \ref{thm_01} is equivalent to the following.
\begin{thm} \label{thm_02}
All zeros of $Z(s)$ lie on the line $Re(s)=1/2$
except for two simple zeros at $s=0, \, 1$.
\end{thm}

\begin{cor} \label{cor_01}
The Riemann zeta function is a factor of the difference of two entire functions whose satisfy
the Riemann hypothesis.
More precisely, we have
\begin{equation} \label{cor_eq1}
s(s-1)\xi(s) \cdot U(s) = V(s) - Z(s),
\end{equation}
where
\begin{equation} \label{cor_eq2}
\aligned
U(s) &= (s+1)\chi(2s-1)-(s-2)\chi(2s), \\
V(s) &= (s-1)(As-A+1)\chi(s+1)\chi(2s) - s(As-1)\chi(s-1)\chi(2s-1).
\endaligned
\end{equation}
\end{cor}
In the proof of Theorem \ref{thm_02}, we use two auxiliary functions
$f(s)$ and $g(s)$ defined by
\begin{equation}\label{06}
\aligned
f(s) & = (s-1)(As-A+1) \cdot \chi(s+1) - (s-2) \cdot \chi(s), \\
g(s) & = f(s) \cdot \chi(2s).
\endaligned
\end{equation}
By definition and \eqref{109}, we have
\begin{equation}
Z(s)=g(s)-g(1-s).
\end{equation}
Here we used the functional equation $\chi(s)=\chi(1-s)$ of the Riemann zeta-function.
\smallskip

Roughly, the proof of Theorem \ref{thm_02} is divided into two steps.
First, we prove that all zeros of $f(s)$
lie in a vertical strip $\sigma_0 < \Re(s) <0$ except for finitely many exceptional zeros (section 2).
Then we obtain a nice product formula of $f(s)$
by a variant of Lemma 3 in~\cite{Su} (Lemma \ref{lem_01}, it is proved in section 4).
Second, by using the product formula of $f(s)$,
we prove that all zeros of $Z(s)$
lie on the line $\Re(s)=1/2$ except for two simple zeros $s=0,1$ (section 3).
In this process, we use the result of Lagarias~\cite{La99} concerning the explicit upper bound
for the difference of the imaginary parts of the zeros of the Riemann zeta function.
In the final section (section 5) we give the proof of Corollary \ref{cor_01}.
\smallskip

Finally, we comment on high rank zeta functions $\xi_{{\Q},r}(s)$
and new zeta functions $\xi_{{\rm SL}(n),{\Q}}(s):=\xi_{\Q}^{G/P}(s)$ attached to $(G,P)=({\rm SL}(n),P_{n-1,1})$.
Roughly, $\xi_{{\Q},r}(s)$ corresponds to $({\rm Res} \to \int)$-ordered construction,
and new zeta functions $\xi_{{\rm SL}(n),{\Q}}(s)$
corresponds to $(\int \to {\rm Res})$-ordered construction.
Here ``$({\rm Res} \to \int)$-ordered'' means that we first taking the residues then taking the integral,
similarly, ``$(\int \to {\rm Res})$-ordered'' means that we first taking the integral then taking the residues.
We have $\xi_{{\Q},2}(s)=\xi_{{\rm SL}(2),{\Q}}(s)$, since we do not need to take a residue.
However, in general, there is a discrepancy between $\xi_{{\Q},r}(s)$ and $\xi_{{\rm SL}(n),{\Q}}(s)$,
because of the obstruction for the exchanging of $\int$ and ${\rm Res}$.
For example, $\xi_{{\Q},3}(s)$ has only two singularities at $s=0,1$,
but $\xi_{{\rm SL}(3),{\Q}}(s)$ has four singularities at $s=0,1/3,2/3,1$.
However, we expect that the distribution of $\xi_{{\rm SL}(n),{\Q}}(s)$ are quite regular
as well as $\xi_{{\Q},r}(s)$.
In fact, we have the Riemann hypothesis of $\xi_{{\rm SL}(2),{\Q}}(s)$
by $\xi_{{\Q},2}(s)=\xi_{{\rm SL}(2),{\Q}}(s)$ and the result of~\cite{LaSu},
and the author proved the Riemann hypothesis of $\xi_{{\rm SL}(3),{\Q}}(s)$ in (\cite{Su}).
Moreover author checked the Riemann hypothesis of $\xi_{{\rm SL}(n),{\Q}}(s)$ for $n=4,5$ by a computational way,
and it suggests that the Riemann hypothesis holds for $n=4,5$.
As observed in~\cite{Su}, the study of $\xi_{{\rm SL}(3),{\Q}}(s)$
gives an information for the Riemann zeta function.
In addition, Corollary \ref{cor_01} gives a relation between
$\xi_{{\rm Sp}(4),{\Q}}(s)$ and the Riemann zeta function.
As these, the study of $\xi_{F}^{G/P}(s)$ is not only interesting itself
but also suggestive for the study of the Dedekind zeta function .
\smallskip

\noindent
{\bf Acknowledgement:} The author would like to thank to Lin Weng
for his contribution of Appendix (\cite{We4}) and for sharing
with me his progress of general zeta functions.
%
%
\section{First step for Theorem \ref{thm_02}}
The aim of this section is to prove the following proposition.
\begin{prop} \label{prop_01}
Let $f(s)$ be the function defined in \eqref{06}.
Then $f(s)$ has the product formula
\begin{equation}
f(s) = f(0) e^{B^\prime s}
\Bigl( 1-\frac{s}{\rho_0} \Bigl)
\Bigl( 1-\frac{s}{\overline{\rho}_0} \Bigl) \cdot \Pi(s)
\quad (B^\prime \geq 0),
\end{equation}
where
\begin{equation}
\Pi(s) = \prod_{{\beta < 1/2}\atop{0 \not= \beta \in {\Bbb R}}} \Bigl( 1 -\frac{s}{\beta} \Bigr)
\prod_{{\beta<1/2}\atop{\gamma >0}} \left[ \Bigl( 1 -\frac{s}{\rho} \Bigr)\Bigl( 1 -\frac{s}{\overline{\rho}} \Bigr) \right]
\quad (\rho=\beta+i\gamma).
\end{equation}
where $\rho_0$ is a complex zero of $f(s)$ with $\Re(\rho_0)>1/2$,
$\beta$ is the real zeros of $f(s)$ and $\rho=\beta+i\gamma$ are other zeros of $f(s)$.
The product converges absolutely on any compact subset in ${\C}$
if we taken the product with the bracket.
\end{prop}
This is proved by checking
that $f(s)$ satisfies all conditions of Lemma \ref{lem_01} in below.
\begin{lem} \label{lem_01}
Let $F(s)$ be an entire function of genus zero or one.
Suppose that
\begin{enumerate}
\item[{\rm (i)}] $F(s)$ is real on the real axis,
\item[{\rm (ii)}] there exists $\sigma_0 >0$ such
that all zeros of $F(s)$ lie in the vertical strip
\begin{equation}\label{le01}
\sigma_0 < \Re(s) < 1/2
\end{equation}
except for finitely many zeros,
\item[{\rm (iii)}]
the zeros of $F(s)$ are finitely many in the right-half plance $\Re(s) \geq 1/2$,
\item[{\rm (iv)}] there exists $C>0$ such that
\begin{equation}\label{le02}
N(T) \leq C T \log T \quad \text{as} \quad T \to \infty,
\end{equation}
where $N(T)$ is the number of zeros of $F(s)$
satisfying $0 \leq \Im(\rho) <T$.
\item[{\rm (v)}] $F(1-\sigma)/F(\sigma) > 0$ for large $\sigma>0$ and
\begin{equation}\label{le03}
F(1-\sigma)/ F(\sigma) \to 0 \quad \text{as} \quad \sigma \to \infty.
\end{equation}
\end{enumerate}
Then $F(s)$ has the product formula
\begin{equation}
F(s)= C s^m e^{B^\prime s}
\prod_{0 \not=\rho \in {\Bbb R}}\Bigl(1-\frac{s}{\rho}\Bigr)
\prod_{\Im(\rho)>0}
\left[ \Bigl(1-\frac{s}{\rho}\Bigr)\Bigl(1-\frac{s}{\bar{\rho}}\Bigr) \right]
\end{equation}
with
\begin{equation}
B^\prime \geq 0,
\end{equation}
The product in the right-hand side converges absolutely on every compact set in $\C$
if we taken the product with the bracket.
\end{lem}
\begin{rem}
This lemma is a variant of Lemma 3 of ~\cite{Su}.
Differ from Lemma 3 of ~\cite{Su}, we allow $F(s)$ to have finitely many zeros in $\Re(s) \geq 1/2$.
The most important part of this lemma is nonnegativity of $B^\prime$.
\end{rem}
We prove Lemma \ref{lem_01} in section 4.
For $f(s)$ in \eqref{06}, (i) is trivial. Under (ii), condition (iv) is easily proved
by the well-known estimate $|\chi(s)| \leq \exp(C|s|\log|s|)$
and Jensen's formula (see~\cite[\S4.1]{Su} for example).
It remains to prove (ii), (iii) and (v) for $f(s)$.

\subsection{Proof of (v)}
First we show that $f(1-\sigma)/f(\sigma)$ is positive for sufficiently large $\sigma>0$.
Using the functional equation of $\chi(s)$, we have
\[
\frac{f(1-\sigma)}{f(\sigma)}
= \frac
{ \sigma \,(A\sigma -1) \cdot \chi(\sigma-1)
+ (\sigma+1) \cdot \chi(\sigma)
}{
(\sigma-1)(A\sigma - A+1) \cdot \chi(\sigma+1)
- (\sigma-2) \cdot \chi(\sigma) }.
\]
Clearly the numerator is positive for large $\sigma>0$.
By ~\cite[pp.109--110]{LaSu}, we have
$|\chi(2\sigma-1)/\chi(2\sigma)|<1$ for $\sigma>1/2$.
Replacing $2\sigma-1$ by $\sigma$
\begin{equation} \label{206}
|\chi(\sigma)/\chi(\sigma+1)|<1 \quad (\sigma>0).
\end{equation}
Hence the denominator is also positive for large $\sigma>0$,
since $A=\pi/3-1>0$.
Now we prove \eqref{le03}. We have
\[
\aligned
\frac{f(1-\sigma)}{f(\sigma)}
& =
\frac{\sigma(A \sigma -1)}{(\sigma-1)(A\sigma-A+1)} \cdot
\frac{\chi(\sigma-1)}{\chi(\sigma+1)} \cdot
\frac{1+g(\sigma)}{1-h(\sigma)} \\
& =
\bigl( 1+O(\sigma^{-1}) \bigr) \cdot
\frac{\chi(\sigma-1)}{\chi(\sigma+1)} \cdot
\frac{1+g(\sigma)}{1-h(\sigma)},
\endaligned
\]
where
\[
g(\sigma)
=
\frac{\sigma+1}{\sigma(A\sigma-1)} \cdot
\frac{\chi(\sigma)}{\chi(\sigma-1)},
\quad
h(\sigma)
=
\frac{\sigma-2}{(\sigma-1)(A\sigma-A+1)} \cdot
\frac{\chi(\sigma)}{\chi(\sigma+1)}.
\]
We have
\[
\aligned
\frac{\chi(\sigma-1)}{\chi(\sigma+1)}
& = \bigl(1+O(\sigma^{-1})\bigr)
\frac{\xi(\sigma-1)}
{\xi(\sigma+1)}
 = \bigl(1+O(\sigma^{-1})\bigr) \cdot \pi \cdot
\frac{\Gamma((\sigma-1)/2)\zeta(\sigma-1)}
{\Gamma((\sigma+1)/2) \zeta(\sigma+1)} \\
& = \bigl(1+O(\sigma^{-1})\bigr) \cdot
\frac{\Gamma((\sigma-1)/2)}
{\Gamma((\sigma+1)/2)} \cdot O(1)
\endaligned
\]
for large $\sigma>0$. Usign the Stirling formula
\[
\Gamma(z)=\sqrt{\frac{2\pi}{z}} \Bigl(\, \frac{z}{e} \,\Bigr)^z \bigl( 1 + O(|z|^{-1})  \bigr)
\quad (|z| \geq 1, ~|{\rm arg}\,z|<\pi-\varepsilon ),
\]
we obtain
\begin{equation} \label{09}
\frac{\chi(\sigma-1)}{\chi(\sigma+1)}
= O(\sigma^{-1}) \quad \text{as} \quad \sigma \to +\infty.
\end{equation}
For $g(\sigma)$, we have
\[
\aligned
g(\sigma)
& =
\frac{\sigma+1}{\sigma(A\sigma-1)} \cdot (1+O(\sigma^{-1})) \cdot \pi^{-1/2}
\frac{\Gamma(\sigma/2)\zeta(\sigma)}
{\Gamma((\sigma-1)/2)\zeta(\sigma-1)} \\
& =
\frac{\sigma+1}{\sigma(A\sigma-1)} \cdot (1+O(\sigma^{-1})) \cdot
\frac{\Gamma(\sigma/2)}
{\Gamma((\sigma-1)/2)} \cdot O(1) \\
& =
\frac{\sigma+1}{\sigma(A\sigma-1)} \cdot (1+O(\sigma^{-1}))
\cdot O(\sqrt{\sigma}).
\endaligned
\]
Here we used the Stirling formula in third equation.
Thus
\begin{equation} \label{10}
g(\sigma) =O(\sigma^{-1/2}) \quad \text{as} \quad \sigma \to +\infty.
\end{equation}
For $h(\sigma)$, by using \eqref{206}, we have
\begin{equation} \label{11}
h(\sigma)=O(\sigma^{-1}) \quad \text{as} \quad \sigma \to +\infty.
\end{equation}
From \eqref{09}, \eqref{10} and \eqref{11}, we obtain
\begin{equation*}
\frac{f(1-\sigma)}{f(\sigma)} =O(\sigma^{-1}) \quad \text{as} \quad \sigma \to +\infty.
\end{equation*}
This shows that $f(s)$ satisfies (v). \hfill $\Box$

\subsection{Proof of (ii) and (iii)}

\begin{lem} \label{lem_02}
The function $f(s)$ in \eqref{06} has no zero in certain left-half plane $\Re(s)<\sigma_1$.
\end{lem}

\noindent
{\bf Proof.}
Assume $\sigma=\Re(s)<0$.
We have
\[
f(s)  = (s-1)(As-A+1)\chi(s+1)
\left[\,
1 -  R(s)
\,\right]
\]
with
\[
R(s) = \frac{s-2}{(s-1)(As-A+1)}
\cdot \frac{\chi(s)}{\chi(s+1)}.
\]
The factor $(s-1)(As-A+1)\chi(s+1)$ has no zero in the left-half plane
$\Re(s)< 1-(1/A) \simeq -20.187$, because of the Euler product of $\chi(s)$.
Using the functional equation of $\chi(s)$,
\[
\aligned
R(s)
& =
\frac{s-2}{(s-1)(As-A+1)}\frac{\chi(1-s)}{\chi(-s)} =
\frac{s-2}{(s+1)(As-A+1)}
\frac{\xi(1-s)}{\xi(-s)} \\
& =
\frac{s-2}{(s+1)(As-A+1)}
\frac{\Gamma((1-s)/2)}{ \sqrt{\pi} \, \Gamma(-s/2)}
\frac{\zeta(1-s)}{\zeta(-s)}.
\endaligned
\]
Therefore
\[
|R(s)|
\leq \frac{1}{\sqrt{\pi}}
\left|
\frac{s-2}{(s+1)(As-A+1)}
\right|
\left|
\frac{\Gamma((1-s)/2)}{\Gamma(-s/2)}
\right|
\zeta(-\sigma)\zeta(1-\sigma).
\]
If $\sigma=\Re(s)<0$, we have $|{\rm arg}(-s/2)| < \pi/2$ and $|{\rm arg}(1-s)/2| < \pi/2$.
Hence we can apply the Stirling formula for $\Re(s)<0$, and then obtain
\[
\left|
\frac{\Gamma((1-s)/2)}{\Gamma(-s/2)}
\right|
= \frac{|s|^{1/2}}{\sqrt{2}}  \cdot \left| 1 -\frac{1}{s} \right|^{-1} \cdot
\frac{1+O(|s|^{-1})}{1+O(|s|^{-1})}=O(|s|^{1/2}).
\]
On the other hand
\[
\zeta(-\sigma)\zeta(1-\sigma) \to 1 \quad (\sigma \to -\infty).
\]
Therefore
\[
|R(s)|
= \frac{1}{\sqrt{2\pi}} \cdot
\frac{|s|^{1/2}}{|(As-A+1)|}(1+o(1))
= \frac{1}{A \sqrt{2\pi|s|}}(1+o(1)),
\]
if $\sigma=\Re(s)<0$, and $|s|$, $|\sigma|$ are both large.
This implies Lemma \ref{lem_02}.
\hfill $\Box$
\smallskip

\begin{lem} \label{lem_03}
The entire function $f(s)$ in \eqref{06} has only finitely many zeros in the right-half plane $\Re(s) > 0$.
In particular, the number of zeros of $f(s)$ in $\Re(s) \geq 1/2$ is finite.
\end{lem}

\noindent
{\bf Proof.}
We have
\begin{equation} \label{14}
f(s)  = (s-1)(As-A+1)\chi(s+1)
\left[\,
1 -  \frac{s-2}{(s-1)(As-A+1)}
\cdot \frac{\chi(s)}{\chi(s+1)}
\,\right].
\end{equation}
By~\cite[pp.109--110]{LaSu},
we have $| \chi(2s-1)/\chi(2s)| < 1$ for any $\Re(s)>1/2$.
Replacing $2s-1$ by $s$, we obtain
\begin{equation}\label{16}
\left| \frac{\chi(s)}{\chi(s+1)} \right|  < 1 \quad (\Re(s) > 0).
\end{equation}
Let $D$ be the region
\[
D:=
\left\{
s \in {\mathbb C} ~\left|~
\Re(s) \geq 0, ~\left| \frac{s-2}{(s-1)(As-A+1)} \right| \geq 1
\right.\right\}.
\]
Then $f(s) \not=0$ if $s \not\in D$ and $\Re(s) \geq 0$,
because of \eqref{14} and \eqref{16}.
The region $D$ is bounded, since
\[
\left| \frac{s-2}{(s-1)(As-A+1)} \right| < 1
\]
for large $|s|$. Hence the number of zeros of $f(s)$ in $\Re(s) \geq 0$ is finite.
\hfill $\Box$
\bigskip

\noindent
{\bf Proof of (ii) and (iii).}
Lemma \ref{lem_02} and Lemma \ref{lem_03} show that (ii) and (iii) hold for $f(s)$.
\hfill $\Box$

\subsection{Proof of Proposition \ref{prop_01}}

By the results in \S 2,1 and \S 2.2, we can apply Lemma \ref{lem_01} to $f(s)$.
On the other hand, we have $f(0)\simeq 1.047 \not=0$.
Hence the proof of Proposition \ref{prop_01} is completed by the following lemma.

\begin{lem} \label{lem_04}
The number of zeros of $f(s)$ in $\Re(s) \geq 1/2$ is just two.
They are non-real zeros and conjugate each other.
The values of them are about $s \simeq 0.927 \pm i \cdot 3.20$.
\end{lem}
\begin{rem}
This lemma is used to simplify the proof of Theorem \ref{thm_02}.
However the explicit values of exceptional zeros are not essential in the proof.
We already know explicitly the region that $f(s)$ possibly have a zero.
From this fact, to prove Theorem \ref{thm_02},
it is sufficient that
we know the explicit number of zeros of $f(s)$ in that region.
\end{rem}
\smallskip

\noindent
{\bf Proof.} The domain $D \cap \{\Re(s) \leq 1/2\}$
is contained in the rectangle $R=[1/2,2] \times [-10,10]$,
where $D$ is the region in the proof of Lemma \ref{lem_03}.
Because of the argument principle, the number of zeros of $f(s)$ in $R$ is given by
\[
\frac{1}{2\pi i} \int_{\partial R} \frac{f^\prime}{f}(s) ds.
\]
In particular, the value of this integral is an integer.
Therefore we can check that the value of this integral is just two
by a computational way (for example Mathematica, Maple, PARI$/$GP, etc.).
Hence we conclude that $f(s)$ has just two zeros in the rectangle $R$.
By another computational way,
we find an approximated value of these two zeros
are $s \simeq 0.927 \pm i \cdot 3.20$.
\hfill $\Box$
%
%
%
\section{Second step for Theorem \ref{thm_02}}
{\bf Proof of  Theorem \ref{thm_02}.}
We have the following three assertions.

\begin{prop} \label{prop_02}
$Z(s)$ has no zero in the right-half plane $\Re(s) \geq 20$.
\end{prop}

\begin{prop} \label{prop_03}
$Z(s)$ has no zero in the region $1/2 < \sigma < 20$, $|t| \geq 22$.
\end{prop}

\begin{prop} \label{prop_04}
$Z(s)$ has only one simple zero $s=1$ in the region $1/2 < \sigma < 20$, $|t| \leq 22$.
\end{prop}
Then, as a consequence of these three results and the functional equation of $Z(s)$,
all zeros of $Z(s)$ lie on the line $\Re(s)=1/2$
except for simple zeros $s=0, \, 1/2, \, 1$. \hfill $\Box$
\bigskip

In the above three assertions, the hardest one is the proof of Proposition \ref{prop_03}.
To prove Proposition \ref{prop_03}, we use the results in the first step
and a result of Lagarias~\cite{La99}.

\subsection{Proof of Proposition \ref{prop_02}}
We have
\begin{equation} \label{3_term}
Z(s)  = (s-1)(As-A+1) \cdot \chi(s+1) \chi(2s) (1 - R_1(s) - R_2(s) - R_3(s)),
\end{equation}
where
\begin{equation}
\aligned
R_1(s) & =  \frac{ (s-2) \cdot \chi(s)}{ (s-1)(As-A+1) \cdot \chi(s+1) }, \\
R_2(s) & =  \frac{ s (As-1) \cdot \chi(s-1) \chi(2s-1) }{ (s-1)(As-A+1) \cdot \chi(s+1) \chi(2s) }, \\
R_3(s) & =  \frac{ (s+1) \cdot \chi(s)  \chi(2s-1) }{ (s-1)(As-A+1) \cdot \chi(s+1) \chi(2s) }.
\endaligned
\end{equation}
For each $R_i(s)~(i=1,2,3)$, we have
\begin{equation}
\aligned
|R_1(s)|
& =  \left| \frac{s-2}{(s+1)(As-A+1)} \frac{ \sqrt{\pi} \, \Gamma(s/2)}{\Gamma((s+1)/2)} \frac{\zeta(s)}{\zeta(s+1)} \right| \\
& \leq  \pi^{1/2} \, \left| \frac{s-2}{(s+1)(As-A+1)} \right|
\left| \frac{ \Gamma(s/2)}{\Gamma((s+1)/2)} \right|
\zeta(\sigma)\zeta(\sigma+1),
\endaligned
\end{equation}
\begin{equation}
\aligned
|&R_2(s)|
 = \pi^{3/2} \, \left| \frac{ (s-1)(s-2) (As-1)}{ s (s+1) (As-A+1)}
\frac{  \Gamma((s-1)/2)\Gamma(s-1/2)}{\Gamma((s+1)/2)\Gamma(s)}
\frac{\zeta(s-1)\zeta(2s-1)}{\zeta(s+1)\zeta(2s)} \right| \\
& \leq 2 \pi^{3/2} \left| \frac{ (s-2) (As-1)}{ s (s+1) (As-A+1)} \right|
\left| \frac{ \Gamma(s-1/2)}{\Gamma(s)}  \right|
\zeta(\sigma+1)\zeta(\sigma-1)\zeta(2\sigma)\zeta(2\sigma-1),
\endaligned
\end{equation}
and
\begin{equation}
\aligned
|R_3(s)|
& = \pi \, \left| \frac{ s-1 }{ s(As-A+1) }
\frac{  \Gamma(s/2) \Gamma(s-1/2) }{ \Gamma((s+1)/2)\Gamma(s) }
\frac{\zeta(s)  \zeta(2s-1)}{\zeta(s+1)\zeta(2s)} \right| \\
& \leq \pi \, \left| \frac{ s-1}{ s(As-A+1)} \right|
\left| \frac{ \Gamma(s-1/2)}{\Gamma(s)}  \right|
\zeta(\sigma+1)\zeta(\sigma)\zeta(2\sigma)\zeta(2\sigma-1).
\endaligned
\end{equation}
Using the Stirling formula, we obtain
\begin{equation}
|R_1(s)| =O(|s|^{-1/2}), \quad |R_2(s)|=(|s|^{-3/2}), \quad |R_3(s)|=(|s|^{-1})
\end{equation}
as $|s|\to \infty$ in the right-half plane.
Therefore $Z(s) \not=0$ for some right-half plane $\Re(s) \geq \sigma_2$.
Using the monotone decreasing property of $\zeta(\sigma)$ as $\sigma \to +\infty$
and the effective version of Stirling's formula (\cite{Ol})
\begin{equation}
\Gamma(s) = \Bigl( \frac{2\pi}{s} \Bigr)^{\frac{1}{2}} \Bigl( \frac{s}{e} \Bigr)^s
\Bigl\{ 1+\Theta \Bigl( \frac{1}{8|s|}\Bigr) \Bigr\} \quad (\Re(s)>1),
\end{equation}
where the notation $f=\Theta(g)$ means $|f| \leq g$,
we have
\begin{equation}
|R_1(s)| \leq 0.5, \quad |R_2(s)| \leq 0.1, \quad |R_3(s)| \leq 0.3
\end{equation}
for $\Re(s) \geq 20$ (in fact, these bounds already hold for $\Re(s) \geq 10$).
These estimates imply $Z(s)\not=0$ for $\Re(s) \geq 20$ by \eqref{3_term},
since $(s-1)(As-A+1) \cdot \chi(s+1) \chi(2s)$ has no zero in the right-half plane $\Re(s) \geq 20$.
\hfill $\Box$

\subsection{Proof of Proposition \ref{prop_03}}
Let $\rho_0=\beta_0+ i \gamma_0~(\gamma_0>0)$ be the zero of $f(s)$ in Lemma \ref{lem_04}.
By Proposition \ref{prop_01} we have
\begin{equation} \label{17}
f(s) = f(0) e^{B^\prime s}
\Bigl( 1-\frac{s}{\rho_0} \Bigl)
\Bigl( 1-\frac{s}{\overline{\rho}_0} \Bigl) \cdot \Pi(s)
\quad (B^\prime \geq 0),
\end{equation}
where
\begin{equation} \label{18}
\Pi(s) = \prod_{{\beta < 1/2}\atop{0 \not= \beta \in {\Bbb R}}} \Bigl( 1 -\frac{s}{\beta} \Bigr)
\prod_{{\rho=\beta+i\gamma}\atop{\beta<1/2,\, \gamma >0}} \left[ \Bigl( 1 -\frac{s}{\rho} \Bigr)\Bigl( 1 -\frac{s}{\overline{\rho}} \Bigr) \right]
\quad (\rho=\beta+i\gamma).
\end{equation}
Having in mind that all zeros of $\Pi(s)$ lie in $\sigma_0 < \Re(s) <1/2$ for some $\sigma_0$
by Lemma \ref{lem_03} and Lemma \ref{lem_04}.
By definition \eqref{06} of $g(s)$,
\begin{equation} \label{28}
Z(s) = g(s) \cdot \Bigl( 1-\frac{g(1-s)}{g(s)} \Bigr) \quad (g(s)=f(s)\cdot \chi(2s)).
\end{equation}
On the other hand,  by \eqref{17},
\begin{equation} \label{29}
\Bigl| \frac{g(1-s)}{g(s)} \Bigr|
= e^{B^\prime (1-2\sigma)}
\cdot \Bigl| \frac{\Pi(1-s)}{\Pi(s)} \Bigr|
\cdot \Bigl| \frac{s-1+\rho_0}{s-\rho_0}
\cdot  \frac{s-1+\overline{\rho}_0}{s-\overline{\rho}_0} \Bigr|
\cdot \Bigl| \frac{\chi(2s-1)}{\chi(2s)} \Bigr|.
\end{equation}
Because $B^\prime \geq 0$ by Lemma \ref{lem_01},
\begin{equation} \label{30}
e^{B^\prime(1-2\sigma)} \leq 1 \quad (\Re(s)>1/2).
\end{equation}
For the ratio of $\Pi(s)$ in \eqref{29}, we have
\begin{equation} \label{31}
\Bigl| \frac{\Pi(1-s)}{\Pi(s)} \Bigr| =
\prod_{{\rho=\beta+i\gamma}\atop{\beta<1/2,\, \gamma >0}}
\left(
\Bigl| \frac{1-s-\overline{\rho}}{s-\rho} \Bigr| \cdot
\Bigl| \frac{1-s-\rho}{s-\overline{\rho}} \Bigr|
\right) < 1 \quad (\Re(s)>1/2),
\end{equation}
by term-by-term argument as in ~\cite{LaSu}
by using $\beta<1/2$ and
\[
\Bigl| \frac{1-s-\overline{\rho}}{s-\rho} \Bigr|^2
= 1 - \frac{(2\sigma-1)(1-2\beta)}{(\sigma-\beta)^2 + (t-\gamma)^2},
\]
where $\rho=\beta+i\gamma$ is a zero of $f(s)$.
It remains to give an estimate for
\begin{equation}
r(s) := \Bigl| \frac{s-1+\rho_0}{s-\rho_0}
\cdot  \frac{s-1+\overline{\rho}_0}{s-\overline{\rho}_0} \Bigr|
\cdot \Bigl| \frac{\chi(2s-1)}{\chi(2s)} \Bigr|.
\end{equation}
To evaluate $r(s)$, we use the following lemma essentially.
\begin{lem}[\cite{La99}] \label{lem_08}
For any real $t$ with $|t| \ge 12$ there
exists at least two distinct zeros $\rho= \beta+ i \gamma$ of $\xi(s)$ such that
$0 < \beta \leq 1/2$ and
\begin{equation} \label{35}
|t - \gamma| \leq 10.1.
\end{equation}
\end{lem}
{\bf Proof.} Suppose $t \geq 25$.
Then there exists at least two distinct zeros $\rho=\beta+i\gamma$ of $\xi(s)$
satisfying $0 < \beta \leq 1/2$ and $|t-\gamma| <10.1$
by applying  Lemma 5 in \cite{Su} to $t+5.1$ and $t-5.1$
(Lemma 5 in \cite{Su} is essentially Lemma 3.5 of \cite{La99}).
For $12 \leq t < 25$, estimate \eqref{35} also holds for two zeros
because $\xi(s)$ has zeros at $s=\pm 14.13, \, \pm 21.02, \, \pm 25.01$.
\hfill $\Box$
\bigskip

\noindent
Using Lemma \ref{lem_08} we show the following.

\begin{lem} \label{lem_09}
Let $\rho_0 = \beta_0 + i \gamma_0 \simeq 0.927 + i \cdot 3.20$ be the zero of $f(s)$ in Lemma \ref{lem_04}.
Let $s = \sigma + it$ with $1/2 < \sigma \leq 20$ and $t \ge 22$.
Then there exists at least two distinct zeros $\rho=\beta+i\gamma$ of $\xi(s)$ such that
$0 < \beta \leq 1/2$, $|t-\gamma| \leq 10.1$,
\begin{equation} \label{35_1}
\left| \frac{s-1+\overline{\rho_0}}{s-\rho_0} \right| \cdot
\left| \frac{2s-1 - ( 1 - \overline{\rho}) }{2s - \rho} \right|
< 1
\end{equation}
and
\begin{equation} \label{35_2}
\left| \frac{s-1+\rho_0}{s-\overline{\rho_0}} \right| \cdot
\left| \frac{2s-1 - ( 1 - \overline{\rho}) }{2s - \rho} \right|
< 1
\end{equation}
\end{lem}
{\bf Proof.}
By squaring \eqref{35_1} and \eqref{35_2} we have
\begin{equation} \label{35_3}
\frac{(\sigma+\beta_0-1)^2 + (t \pm \gamma_0)^2}{(\sigma-\beta_0)^2 + (t \pm \gamma_0)^2}
\cdot \frac{(2\sigma+\beta-2)^2 + (t-\gamma)^2}{(2\sigma-\beta)^2 + (t-\gamma)^2} <1.
\end{equation}
To prove Lemma \ref{lem_09} it is sufficient that \eqref{35_3} holds
for $0 < \beta \leq 1/2$, $|t-\gamma| < 10.1$, $1/2 < \sigma \leq 20$ and $t \geq 22$,
because of Lemma \ref{lem_08}.
To establish \eqref{35_3} in that conditions
it suffices to show that
\begin{equation} \label{35_4}
\frac{(\sigma+\beta_0-1)^2 + (t \pm \gamma_0)^2}{(\sigma-\beta_0)^2 + (t \pm \gamma_0)^2}
\cdot \frac{(2\sigma-\frac{3}{2})^2 + 10^2}{(2\sigma-\frac{1}{2})^2 + 10^2} <1.
\end{equation}
by a similar reason in the later half of \S4.3 in~\cite{Su}.
This inequality is equivalent to
\begin{equation} \label{35_5}
(2\sigma-1) \Bigl( 8(t \pm \gamma_0)^2 - P(\sigma) \Bigr)  > 0,
\end{equation}
where $P(\sigma)=8(4\beta_0-3) \sigma^2 - 8(4\beta_0-3) \sigma - 8 \beta_0^2 + 818 \beta_0 -  409$.
Using the value $\beta_0 \simeq 0.927$ we see that $P(\sigma) < 2580$ for $1/2 < \sigma <20$.
On the other hand, using the value $\gamma_0 \simeq 3.20$ we see that
$8(t \pm \gamma_0)^2 > 2590$ for $t \geq 22$
since $|t \pm \gamma_0| = t \pm \gamma_0>18$ for $t \geq 22$.
Hence \eqref{35_5} hold, and it implies \eqref{35_3}.
\hfill $\Box$
\bigskip

\noindent
Lemma \ref{lem_09} and $\overline{Z(s)}=Z(\overline{s})$ implies
\begin{equation} \label{36}
|r(s)| <1 \quad \text{for} \quad 1/2 < \sigma \leq 20, ~|t| \geq 22
\end{equation}
by taking two distinct zeros of $\xi(s)$ in this region,
since we have the inequality
\begin{equation}
\Bigl|  \frac{2s-1 -(1-\overline{\rho})}{2s -\rho} \Bigr| < 1 \quad (\Re(s)>1/2),
\end{equation}
for other terms in $r(s)$,
where $\rho$ is a zero of $\xi(s)$ ($0<\Re(\rho)<1$).
Estimates \eqref{30}, \eqref{31} and \eqref{36} show that
\begin{equation}
\Bigl| \frac{g(1-s)}{g(s)} \Bigr| < 1
\end{equation}
for $1/2 < \sigma \leq 20$, $|t| \geq 22$.
By \eqref{28} this estimate implies Proposition \ref{prop_03},
because $g(s)$ has no zero in $1/2 < \sigma \leq 20$, $|t| \geq 22$
by Lemma \ref{lem_03} and Lemma \ref{lem_04}.
\hfill $\Box$

\subsection{Proof of Proposition \ref{prop_04}}
Because the region $1/2 < \sigma \leq 20$, $|t| \leq 22$ is finite,
we can check Proposition \ref{prop_04} by using the help of computer
as in the proof of Lemma \ref{lem_04}. \hfill $\Box$

%
\section{Proof of Lemma \ref{lem_01}}
We prove the lemma only if $F(s)$ has genus one,
since if $F(s)$ has genus zero it is easily proved by a way similar to the case of genus one.
The genus one assumption is equivalent to the Hadamard product factorization
\begin{equation} \label{017}
F(s) = e^{A+Bs} s^m
\prod_{\rho} \Bigl( 1 - \frac{s}{\rho} \,\Bigr) \exp(s/\rho) \quad (m \in {\Bbb Z}_{\geq 0})
\end{equation}
converges absolutely and uniformly
on any compact subsets of $\Bbb C$.
That is also equivalent to $\sum_{\rho} |\rho|^{-2} < \infty.$
Assumption (i) implies the symmetry of the set of zeros under the conjugation $\rho \mapsto \overline{\rho}$.
It follows that the set of zeros $\rho= \beta + i \gamma$,
counted with multiplicity,
is partitioned into  blocks $B(\rho)$ comprising
$\{ \rho, \overline{\rho} \}$ if $\gamma >0$
and
$\{\rho\}$ if $\beta \not=0$ and $\gamma=0$.
Each block is labeled with the unique zero
in it having $\gamma \geq 0$.
Using assumption (ii), we show
\begin{equation} \label{018}
F(s) = s^m e^{A+ B^\prime s}
\prod_{B(\rho)}
\left( \prod_{\rho \in B(\rho)} \Bigl(\, 1 - \frac{s}{\rho} \,\Bigr) \right)
\end{equation}
where the outer product on the right-hand side
converges absolutely and uniformly on any compact
subsets of $\Bbb C$. This assertion holds because the block convergence
factors $\exp(c(B(\rho))s)$ are given by $c(B(\rho)) = 2\beta |\rho|^{-2}$
for $\gamma >0$. Assumption (ii) implies $ |\beta -1/2| < \sigma_0$.
Hence
\[
\sum_{B(\rho)} |c(B(\rho))| \leq
  \sum_{0 \not= \rho :\, {\rm real}} |\rho|^{-1}
+ (2\sigma_0 +1) \sum_{\rho} |\rho|^{-2} < \infty.
\]
Thus the convergence factors $\exp(c(B(\rho))s)$ can be pulled out of the product.
Hence we have $\eqref{018}$ with
\begin{equation}\label{le08}
B^\prime = B + \sum_{B(\rho)}c(B(\rho)).
\end{equation}
Using assumption (iii), (iv) and (v) we show
\begin{equation}\label{le09}
B^\prime  \geq 0.
\end{equation}
By $\eqref{le03}$ in assumption (v) we have
\begin{equation}\label{le13}
{\Bbb R} \ni \log \Bigl( \frac{F(1-\sigma)}{F(\sigma)} \Bigr) \to -\infty
\quad \text{as} \quad \sigma \to +\infty.
\end{equation}
Using $\eqref{018}$ we have
\[
\frac{F(1-\sigma)}{F(\sigma)}
 = e^{B^\prime(1-2\sigma)}
\Bigl( \frac{\sigma-1}{\sigma} \Bigr)^m
\prod_{\rho=\beta \in {\Bbb R}}
\frac{\sigma - 1 + \beta} { \sigma - \beta }
\prod_{{\rho=\beta+i\gamma}\atop{\gamma>0}}
\frac{(\sigma - 1 + \beta)^2 + \gamma^2} { (\sigma - \beta)^2 + \gamma^2 }.
\]
Thus
\begin{equation}\label{le14}
\aligned
\log  \Bigl( \frac{F(1-\sigma)}{F(\sigma)} \Bigr)
& = B^\prime(1-2\sigma)
  + m \log \Bigl( 1-  \frac{1}{\sigma} \Bigr)
  + \sum_{\rho=\beta \in {\Bbb R}}
  \log\Bigl(1 -  \frac{1 - 2\beta} { \sigma - \beta } \Bigr)  \\
& \quad
  +\sum_{{\rho=\beta+i\gamma}\atop{\gamma>0}}
  \log\Bigl( 1- \frac{(1-2\beta)(2\sigma-1)} { (\sigma - \beta)^2 + \gamma^2 }  \,\Bigr).
\endaligned
\end{equation}
Note that
\begin{equation} \label{023}
\log\Bigl( 1- \frac{(1-2\beta)(2\sigma-1)} { (\sigma - \beta)^2 + \gamma^2 } \Bigr) <0
\quad \text{for} \quad \sigma>1/2
\end{equation}
if $\beta<1/2$, and
\[
\log \Bigl( 1-  \frac{1}{\sigma} \Bigr),~
\log\Bigl(1 -  \frac{1 - 2\beta} { \sigma - \beta } \Bigr),~
\log\Bigl( 1- \frac{(1-2\beta)(2\sigma-1)} { (\sigma - \beta)^2 + \gamma^2 } \Bigr)
\to 0 \quad \text{as} \quad \sigma \to +\infty
\]
for any fixed $\rho=\beta+i\gamma$.
By assumption (iii), \eqref{023} holds except for finitely many zeros.
Hence if we suppose $B^\prime <0$, $\eqref{le13}$ and $\eqref{le14}$ implies
\begin{equation}\label{le15}
\left\vert
\sum_{{\rho=\beta+i\gamma}\atop{\gamma>0}}
\log\Bigl( 1- \frac{(1-2\beta)(2\sigma-1)} { (\sigma - \beta)^2 + \gamma^2 }  \,\Bigr)
\right\vert \geq 2 |B^\prime| \sigma
\end{equation}
for large $\sigma>1/2$,
because the number of real zeros is also finite by assumption (ii) and (iii).
On the other hand, for large $\sigma >1/2$, we have
\begin{equation*}
\aligned
\left\vert
\sum_{{\rho=\beta+i\gamma}\atop{\gamma>0}}
\log\Bigl( 1- \frac{(1-2\beta)(2\sigma-1)} { (\sigma - \beta)^2 + \gamma^2 }  \,\Bigr)
\right\vert
& \leq \left\vert
\sum_{{\rho=\beta+i\gamma}\atop{\gamma>0}}
\log\Bigl( 1- \frac{(1-2\sigma_0)(2\sigma-1)} { (\sigma - 1/2)^2 + \gamma^2 }  \,\Bigr)
\right\vert \\
& \ll (2\sigma-1) \sum_{{\rho=\beta+i\gamma}\atop{\gamma>0}}
\frac{1} { (\sigma - 1/2)^2 + \gamma^2 }.
\endaligned
\end{equation*}
The sum in the right-hand side can be written as the Stieltjes integral
\[
\int_{\gamma_0}^{\infty} \frac{dN(t)}{(\sigma-1/2)^2 + t^2}.
\]
Using $\eqref{le02}$ in (iv) we have
\[
\int_{\gamma_0}^{\infty} \frac{dN(t)}{(\sigma-1/2)^2 + t^2}
\ll \int_{\gamma_0}^{\infty} \frac{(\log t) \, dt}{(\sigma-1/2)^2 + t^2}
\ll \frac{\log (\sigma + \gamma_0)}{\sigma-1/2}.
\]
Hence we obtain
\begin{equation}\label{le16}
\left\vert
\sum_{{\rho=\beta+i\gamma}\atop{\gamma>0}}
\log\Bigl( 1- \frac{(1-2\beta)(2\sigma-1)} { (\sigma - \beta)^2 + \gamma^2 }  \,\Bigr)
\right\vert \ll \log(\sigma + \gamma_0)
\end{equation}
for large $\sigma>1/2$.
This contradict $\eqref{le15}$.
Thus we find that $\eqref{le09}$ holds.
\hfill $\Box$

\section{Proof of Corollary \ref{cor_01}}
By Theorem \ref{thm_02}, it remains to show that
all zeros of $V(s)$ in \label{cor_eq1} lie on the line $\Re(s)=1/2$.
Taking $v(s)=(s-1)(As-A+1)\chi(s+1)$ and using the functional equation of $\chi(s)$,
we have $V(s) = v(s)\chi(2s)-v(1-s)\chi(2s-1)$.
All zeros of $v(s)$ lie in the strip $1-(1/A)<\Re(s)<0$ except for the simple zero $s=1$.
Then we find that all zeros of $V(s)$ are lie on the line $\Re(s)=1/2$
by a way similar to section 2 replacing $\rho_0$ by $1$.
\hfill $\Box$

\noindent
Masatoshi Suzuki \\
Department of Mathematics \\
Rikkyo University \\
Nishi-Ikebukuro, Toshima-ku \\
Tokyo 171-8501, Japan \\
\texttt{suzuki@rkmath.rikkyo.ac.jp}

%
\appendix%

\vfill
\pagebreak

\section{\!\!\!\!\!\!ppendix: Zeta Functions for $Sp(2n)$}
\author{\bf  Lin WENG\footnote{
{\bf Acknowledgement.} We would like to thank Henry Kim for bringing to our attention the paper of Diehl. \newline
$~\qquad$This is partially supported by JSPS.}
}
\date{Graduate School of Mathematics, Kyushu Univ, Fukuoka, Japan}
\maketitle
\subsection{Introduction}
Associated to a number field $F$ is the genuine high rank zeta function
$\xi_{F,r}(s)$ for every fixed  $r\in\mathbb Z_{>0}$. Being natural generalizations of (completed)
Dedekind zeta functions, these functions satisfy canonical properties for zetas as well.
Namely, they admit meromorphic continuations to the whole complex $s$-plane,
satisfy the functional equation $\xi_{F,r}(1-s)=\xi_{F,r}(s)$ and have only two singularities,
all simple poles, at $s=0,\,1$. Moreover,
it is known that all zeros of $\xi_{F,2}(s)$ lie on the central line $\mathrm {Re}(s)=\frac{1}{2}$.
(We in fact now expect that the Riemann Hypothesis holds for all $\xi_{F,r}(s)$.)

Recall that $\xi_{F,r}(s)$ is defined by
\[
\xi_{F,r}(s):=\Big(|\Delta_F|\Big)^{\frac{rs}{2}}\int_{\mathcal M_{F,r}}\Big(e^{h^0(F,\Lambda)}-1\Big) \cdot
\big(e^{-s}\big)^{\mathrm{deg}(\Lambda)}\,d\mu(\Lambda),\quad \mathrm{Re}(s)>1
\]
where $\Delta_F$ denotes the discriminant of $F$,
$\mathcal M_{F,r}$ the moduli space of semi-stable $\mathcal O_F$-lattices of rank $r$
(here $\mathcal O_F$ denotes the ring of integers),
$h^0(F,\Lambda)$ and $\mathrm{deg}(\Lambda)$ denote the 0-th geo-arithmetic cohomology
and the arithmetic degree of the lattice $\Lambda$,
and $d\mu(\Lambda)$ a certain naturally associated Tamagawa type measure on $\mathcal M_{F,r}$.
(For details, see \cite{W1, W3} for basic theory, and \cite{LS}, \cite{W2},
see also \cite{H}, for the Riemann Hypothesis arguments.)

Algebraic groups associated to $\mathcal O_F$-lattices are general linear group $GL$ and special linear group $SL$.
A natural question then is whether principal lattices associated to other reductive groups $G$
and their associated zeta functions can be introduced and studied.
In this paper we start with symplectic group $Sp$.
In contrasting with a geo-arithmetic method used for high rank zetas,
the one adopted in this paper is rather analytic.
And to avoid further complication, we here only work out the full details for $Sp(4)$ over $\mathbb Q$,
even a general framework for $Sp$ over $F$ is outlined. As a concrete result,
we obtain a precise formula for the zeta function $\xi_{Sp(4),\mathbb Q}(s)$
associated to $Sp(4)$ over $\mathbb Q$, by studying a certain Siegel-Eisenstein period.

The newly obtained zeta $\xi_{Sp(4),\mathbb Q}(s)$  for $Sp(4)$ over $\mathbb Q$ proves to be canonical as well.
For example, $\xi_{Sp(4),\mathbb Q}(s)$ can be meromorphically extended to the whole complex $s$-plane,
satisfies the standard functional equation
\[
\xi_{Sp(4),\mathbb Q}(1-s)=\xi_{Sp(4),\mathbb Q}(s)
\]
and admits only four singularities, all simple poles,
at $s=-1,\,0,\,1,\,2$. Most importantly, $\xi_{Sp(4),\mathbb Q}(s)$
satisfies the RH, a result due to Suzuki (\cite{S2}):
\[
All\ zeros\ of\ \xi_{Sp(4),\mathbb Q}(s)\ lie\ on\ the\ central\ line\ \mathrm{Re}(s)=1/2.
\]

\subsection{Periods for $Sp(2n)$}

\subsubsection{Siegel-Maa\ss ~Eisenstein Series}

Let $G=Sp(2n)$ with $G(\mathbb R)=Sp(2n,\mathbb R)$ the symplectic group of degree $n$ over $\mathbb R$.
Denote by $\frak S:=\frak S_n$  the so-called Siegel upper half-space of size $n$,
and for any $Z\in \frak S$, write $Z=X+\sqrt{-1}Y$ according to its real and imaginary parts,
so that $Y=\mathrm{Im}\,Z>0$ and $Z^t=Z$ is symmetric.
Moreover, for any $M=\begin{pmatrix} A&B\\ C&D\end{pmatrix}\in Sp(2n,\mathbb R)$,
as usual, set $M\langle Z\rangle:=(AZ+B)\cdot(CZ+D)^{-1}$
and write $Y(M):=\mathrm{Im} M\langle Z\rangle$.
This defines a natural transitive action of $Sp(2n,\mathbb R)$ on $\frak S_n$.
Note that the stablizer of $\sqrt{-1}E_n$ is simply $SO(2n,\mathbb R)\cap Sp(2n,\mathbb R)$,
consequently, we obtain  the following well-known identification
$Sp(2n,\mathbb R)/SO(2n,\mathbb R)\cap Sp(2n,\mathbb R)\simeq {\frak S}_n$.

Introduce the Siegel modular group $\Gamma_n:=\big\{\mathrm{diag}(\pm 1,\pm1,\cdots,\pm 1)\big\}\backslash Sp(n,\mathbb Z)$.
Let $\frak P_n:=\Big\{\begin{pmatrix} *&*\\ 0&*\end{pmatrix}\in \Gamma\Big\}$ be a subgroup of $\Gamma$
(associated to  a certain standard maximal parabolic subgroup).

Fix $Z\in Sp(2n,\mathbb R)$, define then the associated Siegel-Maa\ss ~Eisenstein series,
or the same, the Siegel-Epstein zeta function by
\[
E_n(Z;s):=\sum_{\gamma\in\frak P_n\backslash\Gamma_n}\frac{|Y|^{-s}}{\|CZ+D\|^{-2s}}.
\]

\subsubsection{Arthur's Analytic Truncation and Eisenstein Period}

As usual, for sympletic group $G=Sp(2n)$, and for a parabolic group $P$
with Levi decomposition $P=MN$ with $M$ the Levi and $N$ the nilpotent,
denote by $\frak a_0$ (resp. $\frak a_P$) the space of characters associated to the Borel (resp. to $P$).
Denote by $\Delta_0$ the corresponding collection of simple roots.
By definition, an element $T\in\frak a_0$ is said to be {\it sufficiently regular}
and denoted by $T\gg0$ if $\langle \alpha,T\rangle\gg 0$ are large enough for all $\alpha\in\Delta_0$.
Fix such a $T$. Let $\phi:\Gamma\backslash \frak S_n\to \mathbb C$ be a smooth function.
We define Arthur's analytic truncation $\wedge^T\phi$
(for $\phi$ with respect to the parameter $T$)
to be the function on $Sp(2n,\mathbb Z)\backslash Sp(2n,\mathbb R)$ given by
\[
\Big(\wedge^T\phi\Big)(Z)
:=\sum_{P:\mathrm{standard}}(-1)^{\mathrm{rank}(P)}\sum_{\delta\in P(\mathbb Z)\backslash Sp(2n,\mathbb Z)}
\phi_P(\delta g)\cdot\hat\tau_P\big(H_P(\delta g)-T\big),
\]
where $\phi_P$ denotes the constant term of $\phi$ along with the standard parabolic subgroup $P$,
$\hat\tau_P$ is the characteristic function of the so-called positive cone in $\frak a_P$,
and $H_P(Z):=\log_Mm_P(Z)$ is an element in $\frak a_P$.
(For unknown notation, all standard, see e.g., \cite{Ar1,Ar2}, \cite{JLR}, and/or \cite{We-1,W3}.)

Fundamental properties of Arthur's truncation may be summarized in the following:
\vskip 0.30cm
\noindent
{\bf Theorem 2.2.1} (\cite{Ar1,Ar2}, see also \cite{OW})
{\it For  a sufficiently positive $T$ in $\frak a_0$, we have
\begin{description}
\item [(1)] $\wedge^T\phi$ is rapidly decreasing, if $\phi$ is an automorphic form on $G(\mathbb Z)\backslash \frak S_n$,
\item [(2)]  $\wedge^T\circ \wedge^T=\wedge^T$,
\item [(3)] $\wedge^T$ is self-adjoint,
\item [(4)]  $\wedge^T{\bold 1}$ is a characteristic function of a compact subset of $G(\mathbb Z)\backslash \frak S_n$.
\end{description}}

Denote by $\frak F(T)$ the compact subset of $G(\mathbb Z)
\backslash \frak S_n$ whose characteristic function is given by
$\wedge^T\bold 1$ in (4).
\vskip 0.30cm
\noindent
{\bf Corollary 2.2.1} (\cite{W1,W3}, see also \cite{KW})
{\it Let $T\gg0$ be a fixed element in $\frak a_0$.
For an automorphic form $\phi$ on $G(\mathbb Z)\backslash \frak S_n$,
\[
\int_{G(\mathbb Z)\backslash \frak S_n}\wedge^T\phi(g)\,
dg=\int_{\frak F(T)}\phi(g)\,dg.
\]}

\noindent
{\it Proof.}
By (1), $\int_{G(\mathbb Z)\backslash \frak S_n}\wedge^T\phi(g)\,dg$ is well-defined. Moreover,
\[
\begin{aligned}
&\int_{G(\mathbb Z)\backslash \frak S_n} \wedge^T\phi(g)\,dg\\
=&\int_{G(\mathbb Z)\backslash \frak S_n}\bold 1(g)\cdot\Big(\wedge^T\circ
\wedge^T\Big)\phi(g)\,dg\qquad(\text{by\ (2)\ above})\\
=&\int_{G(\mathbb Z)\backslash \frak S_n}\wedge^T\bold 1(g)\cdot\wedge^T\phi(g)
\,dg\\
&\quad(\text{by\ (3)\ above\ since}\ \wedge^T\phi(g)\
\text{is\ rapidly\ decreasing})\\
=&\int_{G(\mathbb Z)\backslash \frak S_n}\Big(\wedge^T
\circ\wedge^T\Big)\bold 1(g)\cdot\phi(g)\,dg\\
&\quad(\text{by\ (3)\ again\ since}\ \phi\ \text{is\ of\ moderate\
growth\ and}\ \wedge^T\bold 1\ \text{is\ compactly\ supported})\\
=&\int_{G(\mathbb Z)\backslash \frak S_n}\wedge^T\bold 1(g)\cdot\phi(g)\,dg
\qquad(\text{by\ (2)\ again})\\
=&\int_{\frak F(T)}\phi(g)\,dg
\qquad(\text{by\ (4)}).
\end{aligned}
\]

\subsubsection{Siegel-Maa\ss-period: an analog of high rank zeta}

Motivated by our study on the high rank zeta associated to $SL(n)$ in (\cite{W1,W3}),
we define the {\it Siegel-Maa\ss-period for $Sp(n)$ over $\mathbb Q$},
a special kind of the so-called Eisenstein period, by
\[
Z_{Sp(2n),\mathbb Q}^T(s):=\int_{\Gamma\backslash \frak S_n}\wedge^T E_n(Z;s)\,d\mu(Z).
\]
This is then a function on $s$ depending also on the parameter $T$.
By Corollary 2.2.1,
\[
Z_{Sp(2n),\mathbb Q}^T(s)=\int_{\frak F(T)} E_n(Z;s)\,d\mu(Z).
\]
Thus the study of $Z_{Sp(2n),\mathbb Q}^T(s)$ may be carried out from that for the Siegel-Maa\ss ~series $E_n(Z;s)$.

\subsubsection{Siegel-Eisenstein Series}
In general, it is {\it very difficult}, most of the time, {\it quite impossible,} to calculate Eisenstein periods.
However, if the original automorphic form from which  the Eisenstein series in use is defined is {\it cuspidal},
then  an advanced version of Rankin-Selberg method can be applied to evaluate them. Motivated by this,
we in this subsection explain a method due to Diehl to realize the Siegel-Maa\ss ~series,
which may be viewed as an Eisenstein series associated to the constant function one on a certain {\it maximal parabolic},
as the residue of the so-called Siegel-Eisenstein series associated to the constant function one on the {\it Borel}.

As usual, corresponding to the partition $n=r+1+1+\cdots+1$,
introduce the standard parabolic subgroup
$\frak P_r:=\Big\{\begin{pmatrix} A&*\\ 0&B\end{pmatrix}\in\Gamma\Big\}$,
where $A=\begin{pmatrix}H^t&&&\\ &1&0&\\ &*&\ddots&\\ &&&1\end{pmatrix}$,
$B=\begin{pmatrix}H^{-1}&&&\\ &1&*&\\ &0&\ddots&\\ &&&1\end{pmatrix}$ with $H=H^{(r)}, |H|=1$.
Accordingly, define the associated Siegel-Eisenstein series by
\[
E_r\Big(Z;s_r,\dots,s_n\Big):=\sum_{\gamma\in\frak P_r\backslash \Gamma}\prod_{v=r}^n|Y(\gamma)_v|^{-s_v}.
\]
Here for a matrix $A=(a_{ij})_{i,j=1}^n$,
denote by $A_v$ the matrix $A_v:=(a_{ij})_{i,j=1}^v$, $1\leq v\leq n$.

It is known that such Siegel-Eiesnetsin series are naturally related to the Siegel zeta functions
associated to the standard parabolic subgroup $\frak Q_r$ of $SL(n)$.
More precisely, let $\frak R:=\{\mathrm{diag}(\pm 1,\dots,\pm 1)\}\backslash SL(n,\mathbb Z)$
and $\frak Q_r$ the standard parabolic subgroup associated to the partition $n=r+1+1+\cdots+1$
consisting of matrices in $SL(n,\mathbb Z)$ of the form
$\begin{pmatrix}H&&&\\ &1&*&\\ &0&\ddots&\\ &&&1\end{pmatrix}$ with $H=H^{(r)}, |H|=1$.
Define the associated Siegel zeta functions by
\[
\xi_r^*\Big(Y;s_r,\dots,s_{n-1}\Big):=\sum_{N\in \frak Q_r\backslash \frak R}\prod_{v=r}^{n-1}|Y[N]_v|^{-s_v}
\]
for all $1\leq r\leq n-1$. Then,  we have
\vskip 0.30cm
\noindent
{\bf Lemma 2.4.1} (\cite{D}) {\it With the same notation as above,}

\noindent
(i) $\displaystyle{E_r\Big(Z;s_r,\dots, s_n\Big)=\sum_{\gamma\in \frak B\backslash\Gamma}|Y(\gamma)|^{-s_n}\cdot
\xi_r^*\Big(Y(\gamma);s_r,\dots,s_{n-1}\Big);}$

\noindent
(ii) {\it There exists a constant $c$ depending only on $r$ such that}
\[
\mathrm{Res}_{s_r=\frac{r+1}{2}}\Big(\xi_r^*\Big(Y;s_r,\dots,s_{n-1}\Big)\Big)
= c_r\cdot\xi_{r+1}^*\Big(Y;s_{r+1}+\frac{r}{2},s_{r+2},\dots,s_{n-1}\Big).
\]

Consequently, we have, up to constant factors,
\[
\begin{aligned}
&\mathrm{Res}_{s_{r+1}+\frac{r}{2}
=\frac{(r+1)+1}{2}}\mathrm{Res}_{s_r=\frac{r+1}{2}}\xi_r^*(Y;s_r,\dots,s_{n-1})\\
=&
\mathrm{Res}_{s_{r+1}+\frac{r}{2}=\frac{(r+1)+1}{2}}\xi_{r+1}^*(Y;s_{r+1}+\frac{r}{2},\dots,s_{n-1})\\
=&
\xi_{r+2}^*(Y;s_{r+2}+\frac{r+1}{2},\dots,s_{n-1}).
\end{aligned}
\]
Thus, by taking $r=1$ and repeating this process, we obtain the following
\[
\mathrm{Res}_{s_{n-1}=1}\cdots \mathrm{Res}_{s_{2}=1}\mathrm{Res}_{s_{1}=1}\Big(\xi_1^*(Y;s_1,s_2,\dots, s_{n-1})\Big)
= |Y|^{-\frac{n-1}{2}}
\]
up to a constant factor.
In particular, we get, the following
\vskip 0.30cm
\noindent
{\bf Lemma 2.4.2.} {\it Up to a constant factor,}
\[
\mathrm{Res}_{s_{n-1}=1}\cdots \mathrm{Res}_{s_{2}=1}\mathrm{Res}_{s_{1}=1}\Big(E_r\Big(Z;s_r,\dots, s_n\Big)\Big)
= E_n\Big(Z;s_n+\frac{n-1}{2}\Big).
\]

\noindent
{\it Proof.} Indeed, up to constant factors,
\[
\begin{aligned}
&\mathrm{Res}_{s_{n-1}=1}\cdots \mathrm{Res}_{s_{2}=1}\mathrm{Res}_{s_{1}=1}E_r(Z;s_r,\dots, s_n)\\
=&\sum_{\gamma\in \frak B\backslash\Gamma}|Y(\gamma)|^{-s_n}\cdot
\mathrm{Res}_{s_{n-1}=1}\cdots \mathrm{Res}_{s_{2}=1}\mathrm{Res}_{s_{1}=1}\xi_r^*(Y(\gamma);s_r,\dots,s_{n-1})\\
=&\sum_{\gamma\in \frak B\backslash\Gamma}|Y(\gamma)|^{-s_n}\cdot
|Y(\gamma)|^{-\frac{n-1}{2}}=E_n(Z;s_n+\frac{n-1}{2}).
\end{aligned}
\]

\subsubsection{Relation with Langlands' Eisenstein Series}
To facilitate further discussions, we next write classical Siegel Eisenstein series in terms of Langlands' language.

Let $\lambda=(z_1,z_2,\dots,z_n)\in\frak a_0$, then for $Z=X+\sqrt{-1}Y\in \frak S$, set
\[
{\bf a}^{\lambda}(Z)=\prod_{v=1}^na_v^{-z_v}\qquad \mathrm{with}\qquad a_v=|Y_v|/|Y_{v-1}|.
\]
and the Langlands-Eisenstein series associated to the constant function one on the Borel is defined by
\[
E({\bf 1};\lambda;Z):=\sum_{\delta\in B(\mathbb Z)\backslash \Gamma}{\bf a}^{\lambda}(\delta Z).
\]
As such, then the so-called power function
\[
{\bf p}_{-\bf s}(Y):=\prod_{\mu=1}^n|Y_\mu|^{-s_\mu}
\]
is given by
\[
\begin{aligned}
\prod_{\mu=1}^n|Y_\mu|^{-s_\mu}=&{\bf p}_{-\bf s}(Y)={\bf a}^{\lambda}(Y)
=\prod_{v=1}^na_v^{-z_v}\\
=&|Y_1|^{-z_1+z_2}|Y_2|^{-z_2+2_3}\cdots|Y_{n-1}|^{-z_{n-1}+z_n}|Y_n|^{-z_n}.
\end{aligned}
\]
That then gives the following relations among $z_i$'s and $s_j$'s:
\[
\begin{cases}
s_1=z_1-z_2,&\\
s_2=z_2-z_3,&\\
\dots&\\
s_{n-1}=z_{n-1}-z_n,&\\
s_n=z_n.&
\end{cases}
\]
Consequently, by Lemma 2.4.2, we have  the following
\vskip 0.30cm
\noindent
{\bf Lemma 2.5.1} {\it With the same notation as above,}

\noindent
(i) $\displaystyle{E({\bf 1};z_1,z_2,\dots,z_n;Z)=E_1(Z;s_1,s_2,\dots,s_n)};$

\noindent
(ii) {\it Up to a suitable constant factor,
\[
\begin{aligned}
&E_n\Big(Z,z_n+\frac{n-1}{2}\Big)\\
=&\mathrm{Res}_{z_{n-1}-z_n=1}\cdots \mathrm{Res}_{z_{2}-z_3=1}\mathrm{Res}_{z_{1}-z_2=1}
\Big(E\Big({\bf 1};z_1,z_2,\dots, z_n;Z\Big)\Big).
\end{aligned}
\]}
{\it In particular, in the case $n=2$, i.e, for $Sp(4)$, we have}
\[
\mathrm{Res}_{z_{1}-z_2=1}\Big(E\Big({\bf 1};z_1,z_2;Z \Big)\Big)=E_2\Big(Z,z_2+\frac{1}{2}\Big).
\]

\subsubsection{Advanced Rankin-Selberg Method}

The advantage of using $E({\bf 1};z_1,z_2,\dots,z_n;Z)$ instead of directly using $E_n(Z,s)$
is that the Eisenstein periods for $E({\bf 1};z_1,z_2,\dots,z_n;Z)$ can be evaluated.
Indeed, if $\varphi$ is an automorphic form of  $P'=M'N'$-level,
where $P'$ is a standard parabolic subgroup of a reductive group $G$,
then we can form the associated Eisenstein series
\[
E(\varphi;\lambda;g):=\sum_{\delta\in P'(\mathbb Z)\backslash G(\mathbb Z)}
m_{P'}(\delta g)^{\lambda+\rho}\cdot\varphi(\delta g),\qquad\lambda\in{\mathcal C}_{P'}.
\]
\vskip 0.30cm
\noindent
{\bf Theorem 2.6.1} (\cite[Corollary 17]{JLR})
{\it
Let $P=MN$ be a minimal parabolic subgroup and let $\varphi$ be a $P$-level cusp form.
Let $E(\varphi; \lambda;g)$ be the Eisenstein series associated to $\varphi$.
Then
\[
\int_{G(\mathbb Z)\backslash G(\mathbb R)} \wedge^T E(g,\varphi,\lambda)\, dg
\]
is equal to
\[
v\sum_{w\in W} \frac {e^{\langle w\lambda-\rho, T\rangle}}{\prod_{\alpha\in \Delta_0}
\langle w\lambda-\rho, \alpha^{\vee}\rangle}
\int_{M(\mathbb Z)\backslash M(\mathbb R)^1 \times K}
M(w,\lambda)\varphi(mk)\, dmdk,
\]
where $v=\mathrm{Vol}\Big(\big\{\sum_{\alpha\in \Delta_0} a_{\alpha} \alpha^{\vee}: 0\leq a_{\alpha}<1\big\}\Big)$,
$W$ denotes the Weyl group and $M(w,\lambda)$ denotes the so-called intertwining operator.}
\vskip 0.30cm
In particular, we have the following
\vskip 0.30cm
\noindent
{\bf Corollary 2.6.1} {\it With the same notation as above, up to a constant factor,}
\[
\begin{aligned}
&\int_{Sp(2n,\mathbb Z)\backslash \frak S}\Big(\wedge^T E\Big({\bf 1};z_1,z_2,\dots, z_n;M\Big)\Big)\, d\mu(M)\\
=&\sum_{w\in W} \frac {e^{\langle w\lambda-\rho, T\rangle}}{\prod_{\alpha\in \Delta_0}
\langle w\lambda-\rho, \alpha^{\vee}\rangle} \cdot \prod_{\alpha>0, w\alpha<0}
\frac {\xi\big(\langle\lambda, \alpha^{\vee}\rangle\big)}{\xi\big(\langle \lambda, \alpha^{\vee}\rangle+1\big)}.
\end{aligned}
\]

\noindent
{\it Proof.} In fact, by the Gindikin-Karpelevich formula,
we have
\[
M(w,\lambda)=\prod_{\alpha>0, w\alpha<0}
\frac {\xi\big(\langle\lambda, \alpha^{\vee}\rangle\big)}{\xi\big(\langle \lambda, \alpha^{\vee}\rangle+1\big)}.
\]
Here $\xi(s)$ is the completed Riemann zeta function with the usual $\Gamma$-factors, namely,
$\xi(s)=\pi^{-\frac s2}\Gamma(\frac s2)\zeta(s)$, where $\zeta(s)$ is the Riemann zeta function.

\subsubsection{Periods for $Sp(2n)$ over $\mathbb Q$}

Recall that the Siegel-Maa\ss-period $Z_{Sp(2n),\mathbb Q}^T(s)$,
an analog of high rank zeta for $Sp(2n)$, is equal to
\[
\int_{\frak F(T)}E_n(Z;s)\,d\mu(Z).
\]
Thus to evaluate the Siegel-Maa\ss-period $Z_{Sp(2n),\mathbb Q}^T(s)$,
it suffices to evaluate the integration
\[
\int_{\frak F(T)}\mathrm{Res}_{z_{n-1}-z_n=1}\cdots \mathrm{Res}_{z_{2}-z_3=1}\mathrm{Res}_{z_{1}-z_2=1}
\Big(E\Big({\bf 1};z_1,z_2,\dots, z_n;Z\Big)\Big)\,d\mu(Z).
\]
Consequently, if we were able to freely make an interchange between

\noindent
(i) the operation of taking integration $\int_{\frak F(T)}$ and

\noindent
(ii) the operation of taking residues $\mathrm{Res}_{z_{n-1}-z_n=1}\cdots \mathrm{Res}_{z_{2}-z_3=1}\mathrm{Res}_{z_{1}-z_2=1}$,

\noindent
it would be sufficient for us to evaluate
\[
\mathrm{Res}_{z_{n-1}-z_n=1} \cdots \mathrm{Res}_{z_{2}-z_3=1}\mathrm{Res}_{z_{1}-z_2=1}
\Big(\int_{\frak F(T)}E({\bf 1};z_1,z_2,\dots, z_n;Z)\,d\mu(Z)\Big),
\]
or better, to evaluate the expression
\[
\begin{aligned}
\mathrm{Res}_{z_{n-1}-z_n=1}&\cdots \mathrm{Res}_{z_{2}-z_3=1}\mathrm{Res}_{z_{1}-z_2=1}\\
&\Bigg(\sum_{w\in W} \frac {e^{\langle w\lambda-\rho, T\rangle}}{\prod_{\alpha\in \Delta_0}
\langle w\lambda-\rho, \alpha^{\vee}\rangle} \cdot \prod_{\alpha>0, w\alpha<0}
\frac {\xi\big(\langle\lambda, \alpha^{\vee}\rangle\big)}{\xi\big(\langle \lambda, \alpha^{\vee}\rangle+1\big)}\Bigg)
\end{aligned}
\]
since by Corollary 2.2.1,
\[
\begin{aligned}
\int_{\frak F(T)}&E({\bf 1};z_1,z_2,\dots, z_n;Z)\,d\mu(Z)\\
=&\int_{Sp(2n,\mathbb Z)\backslash \frak S_n}\Lambda^TE({\bf 1};z_1,z_2,\dots, z_n;Z)\,d\mu(Z).
\end{aligned}
\]

Unfortunately, this interchange of orders of two operations is not allowed in general.
As examples, one can observe this by working on $SL(n)$ and by comparing the poles for the resulting expressions.
\vskip 0.30cm
On the other hand, even with the existence of such discrepancies, the function
\[
\begin{aligned}
\mathrm{Res}_{z_{n-1}-z_n=1}&\cdots \mathrm{Res}_{z_{2}-z_3=1}\mathrm{Res}_{z_{1}-z_2=1}\\
&\Bigg(\sum_{w\in W} \frac {e^{\langle w\lambda-\rho, T\rangle}}{\prod_{\alpha\in \Delta_0}
\langle w\lambda-\rho, \alpha^{\vee}\rangle}  \cdot \prod_{\alpha>0, w\alpha<0}
\frac {\xi\big(\langle\lambda, \alpha^{\vee}\rangle\big)}{\xi\big(\langle \lambda, \alpha^{\vee}\rangle+1\big)}\Bigg)
\end{aligned}
\]
proves to be extremely {\it natural and nice}.
To see this, and to make the discussion simpler,
let now concentrate to the case when $T=0$.
\vskip 0.30cm
\noindent
{\bf Definition.} (1) The period $\omega_{Sp(2n),\mathbb Q}(\lambda)$ associated to $Sp(2n)$ over $\mathbb Q$ is defined by
\[
\omega_{Sp(2n),\mathbb Q}(\lambda):=\sum_{w\in W} \Big(\frac {1}{\prod_{\alpha\in \Delta_0}
\langle w\lambda-\rho, \alpha^{\vee}\rangle}  \cdot \prod_{\alpha>0, w\alpha<0}
\frac {\xi\big(\langle\lambda, \alpha^{\vee}\rangle\big)}{\xi\big(\langle \lambda, \alpha^{\vee}\rangle+1\big)}\Big);
\]

\noindent
(2) The period $Z_{\mathbb Q}^{Sp(2n)}(z_n)$ associated to $Sp(2n)$ over $\mathbb Q$ is defined by
\[
\begin{aligned}
Z_{\mathbb Q}^{Sp(2n)}(z_n):=&\mathrm{Res}_{z_{n-1}-z_n=1}\cdots \mathrm{Res}_{z_{2}-z_3=1}\mathrm{Res}_{z_{1}-z_2=1}\\
&\Bigg(\sum_{w\in W} \Big(\frac {1}{\prod_{\alpha\in \Delta_0}
\langle w\lambda-\rho, \alpha^{\vee}\rangle}  \cdot \prod_{\alpha>0, w\alpha<0}
\frac {\xi\big(\langle\lambda, \alpha^{\vee}\rangle\big)}{\xi\big(\langle \lambda, \alpha^{\vee}\rangle+1\big)}\Big)\Bigg).
\end{aligned}
\]

\noindent
{\it Remarks.} \\
(1) Recall that in the original discussions in subsections A.2.2, A.2.3 and A.2.6,
$T$ was assumed to be sufficiently positive, for the reason of Corollary 2.2.1.
However, with the use of the above concrete expression,
$T$ can be chosen to be any element in $\frak a_0$.
\smallskip

\noindent
(2) For high rank zetas,  i.e., in the case where the corresponding algebraic group is $SL(n)$,
the corresponding $\frak F(T)$ makes sense also for $T=0$.
Indeed, $\frak F(0)$ coincides with the moduli space of semi-stable lattices of volume one and rank $n$.
Consequently, the corresponding period after putting $T=0$
gives essentially the high rank zeta $\xi_{\mathbb Q,r}(s)$ there.
So to obtain an analogue of high rank zeta for $Sp(2n)$ geo-arithmetically,
we need to understand $\frak F(0)$ geo-arithmetically.

\subsection{Zetas for $Sp(2n)$}
\subsubsection{An Zeta for $Sp(4)$}
Now we focus on the case $G=Sp(4)$.
Then $\Delta_0=\{e_1-e_2, 2e_2\}$.
There are 4 positive roots $\{e_1\pm e_2, 2e_1, 2e_2\}$ and the Weyl group consists of 8 Weyl elements
\[
\{1, (12), c_1, c_2, (12)c_1, (12)c_2, (12)c_1c_2, c_1c_2\},
\]
where $(12)$ is the lower indices change and $c_i$'s are sign changes.
Also $\rho=2e_1+e_2$.
Consequently, we have the following table for the Weyl action on positive roots:

\[
\begin{matrix}
&e_1-e_2&2e_2&|&e_1+e_2&2e_1\\
1&e_1-e_2&2e_2&|&e_1+e_2&2e_1\\
(12)&e_2-e_1&2e_1&|&e_1+e_2&2e_2\\
c_1&-e_1-e_2&2e_2&|&-e_1+e_2&-2e_1\\
c_2&e_1+e_2&-2e_2&|&e_1-e_2&2e_1\\
(12)c_1&-e_1-e_2&2e_1&|&e_1-e_2&2e_2\\
(12)c_2&e_1+e_2&-2e_1&|&e_2-e_1&2e_2\\
(12)c_1c_2&e_1-e_2&-2e_1&|&-e_1-e_2&-2e_2\\
c_1c_2&e_2-e_1&-2e_2&|&-e_1-e_2&2e_1
\end{matrix}
\]
So to calculate the part
$\frac{1}{\langle\lambda, w^{-1}\alpha_1\rangle-1}\cdot\frac{1}{(\langle\lambda, w^{-1}\alpha_2\rangle-1}$,
we need to use the left half
\[
\begin{matrix}
&e_1-e_2&2e_2\\
1=1^{-1}&e_1-e_2&2e_2\\
(12)=(12)^{-1}&e_2-e_1&2e_1\\
c_1=c_1^{-1}&-e_1-e_2&2e_2\\
c_2=c_2^{-1}&e_1+e_2&-2e_2\\
(12)c_1=\big((12)c_2)\big)^{-1}&-e_1-e_2&2e_1\\
(12)c_2=\big((12)c_1\big)^{-1}&e_1+e_2&-2e_1\\
(12)c_1c_2&e_1-e_2&-2e_1\\
c_1c_2&e_2-e_1&-2e_2
\end{matrix}
\]
and to calculate the zeta part
$\prod_{\alpha>0,w\alpha<0}
\frac {\xi\big(\langle\lambda, \alpha^{\vee}\rangle\big)}{\xi\big(\langle \lambda, \alpha^{\vee}\rangle+1\big)}$,
we need to use the following table:

\[
\begin{matrix}
&e_1-e_2&2e_2&|&e_1+e_2&2e_1\\
1&\times&\times&|&\times&\times\\
(12)&0&\times&|&\times&\times\\
c_1&0&\times&|&0&0\\
c_2&\times&0&|&\times&\times\\
(12)c_1&0&\times&|&\times&0\\
(12)c_2&\times&0&|&0&\times\\
(12)c_1c_2&\times&0&|&0&0\\
c_1c_2&0&0&|&0&0
\end{matrix}
\]
where $\lq\times$' means the corresponding positive root will not contribute
and $\lq 0$' means the corresponding positive root will contribute as it changes to a negative root
under the corresponding Weyl action.
From all this, by a routine but direct calculation, which we decide to omit, we obtain the following table
\[
\begin{matrix}
&\frac{1}{\langle\lambda, w^{-1}\alpha_1\rangle-1}\cdot\frac{1}{\langle\lambda, w^{-1}\alpha_2\rangle-1}&|&\prod_{\alpha>0,w\alpha<0} \frac {\xi\big(\langle\lambda, \alpha^{\vee}\rangle\big)}{\xi\big(\langle \lambda, \alpha^{\vee}\rangle+1\big)}\\
1&\frac{1}{(z_1-z_2-1)}\frac{1}{(z_2-1)}&|&1\\
(12)&\frac{1}{(z_2-z_1-1)}\frac{1}{(z_1-1)}&|&\frac{\xi(z_1-z_2)}{\xi(z_1-z_2+1)}\\
c_1&\frac{1}{(-z_1-z_2-1)}\frac{1}{(z_2-1)}&|&\frac{\xi(z_1-z_2)}{\xi(z_1-z_2+1)}\cdot\frac{\xi(z_1+z_2)}{\xi(z_1+z_2+1)}\cdot\frac{\xi(z_1)}{\xi(z_1+1)}\\
c_2&\frac{1}{(z_1+z_2-1)}\frac{1}{(-z_2-1)}&|&\frac{\xi(z_2)}{\xi(z_2+1)}\\
(12)c_1&\frac{1}{(z_1+z_2-1)}\frac{1}{(-z_1-1)}&|&\frac{\xi(z_1-z_2)}{\xi(z_1-z_2+1)}\cdot\frac{\xi(z_1)}{\xi(z_1+1)}\\
(12)c_2&\frac{1}{(-z_1-z_2-1)}\frac{1}{(z_1-1)}&|&\frac{\xi(z_1+z_2)}{\xi(z_1+z_2+1)}\cdot\frac{\xi(z_2)}{\xi(z_2+1)}\\
(12)c_1c_2&\frac{1}{(z_1-z_2-1)}\frac{1}{(-z_1-1)}&|&\frac{\xi(z_2)}{\xi(z_2+1)}\cdot\frac{\xi(z_1+z_2)}{\xi(z_1+z_2+1)}\cdot\frac{\xi(z_1)}{\xi(z_1+1)}\\
c_1c_2&\frac{1}{(-z_1+z_2-1)}\frac{1}{(-z_2-1)}&|&\frac{\xi(z_1-z_2)}{\xi(z_1-z_2+1)}\cdot\frac{\xi(z_2)}{\xi(z_2+1)}\cdot\frac{\xi(z_1+z_2)}{\xi(z_1+z_2+1)}\cdot\frac{\xi(z_1)}{\xi(z_1+1)}.
\end{matrix}
\]
As such, by taking the residue along the singular line $z_1-z_2=1$ and setting $z_1=b+1, z_2=b$, we get, for the product
\[
\frac{1}{\langle\lambda,
w^{-1}\alpha_1\rangle-1}\cdot\frac{1}{\langle\lambda, w^{-1}\alpha_2\rangle-1}\cdot\prod_{\alpha>0,w\alpha<0}
\frac {\xi\big(\langle\lambda, \alpha^{\vee}\rangle\big)}{\xi\big(\langle \lambda, \alpha^{\vee}\rangle+1\big)},
\]
the following table of contributions
\[
\begin{matrix}
1&\frac{1}{b-1}\\
(12)&\frac{1}{-2}\frac{1}{b}\frac{1}{\xi(2)}\\
c_1&\frac{1}{-2b-2}\frac{1}{b-1}\frac{1}{\xi(2)}\frac{\xi(2b+1)}{2b+2)}\frac{\xi(b+1)}{\xi(b+2)}\\
c_2&0\\
(12)c_1&\frac{1}{2b}\frac{1}{-b-2}\frac{1}{\xi(2)}\frac{\xi(b+1)}{\xi(b+2)}\\
(12)c_2&0\\
(12)c_1c_2&\frac{1}{-b-2}\frac{\xi(b)}{\xi(b+1)}\frac{\xi(2b+1)}{\xi(2b+2)}\frac{\xi(b+1)}{\xi(b+2)}\\
c_1c_2&\frac{1}{-2}\frac{1}{-b-1}\frac{1}{\xi(2)}
\frac{\xi(b)}{\xi(b+1)}\frac{\xi(2b+1)}{\xi(2b+2)}\frac{\xi(b+1)}{\xi(b+2)}
\end{matrix}
\]
This then leads to the following explicit expression of the period $Z_{\mathbb Q}^{Sp(4)}$:
\[
\begin{aligned}
Z_{\mathbb Q}^{Sp(4)}(b)=&\frac{1}{b-1}+\frac{1}{-2}\frac{1}{b}\frac{1}{\xi(2)}\\
&+\frac{1}{-2b-2}\frac{1}{b-1}\frac{1}{\xi(2)}\frac{\xi(2b+1)}{2b+2)}\frac{\xi(b+1)}{\xi(b+2)}\\
&+\frac{1}{2b}\frac{1}{-b-2}\frac{1}{\xi(2)}\frac{\xi(b+1)}{\xi(b+2)}\\
&+\frac{1}{-b-2}\frac{\xi(b)}{\xi(b+1)}\frac{\xi(2b+1)}{\xi(2b+2)}\frac{\xi(b+1)}{\xi(b+2)}\\
&+\frac{1}{-2}\frac{1}{-b-1}\frac{1}{\xi(2)}
\frac{\xi(b)}{\xi(b+1)}\frac{\xi(2b+1)}{\xi(2b+2)}\frac{\xi(b+1)}{\xi(b+2)}.
\end{aligned}
\]

Now multiplying the period $Z_{\mathbb Q}^{Sp(4)}(b)$ with the factor $\xi(b+2)\cdot\xi(2b+2)$,
for the purpose of clearing up the xi function factors appeared in the denominators,
and with the factor $\xi(2)$ for the purpose of clearing up the xi special values appeared in the denominators,
we then obtain the following function
\[
\xi_{Sp(4),\mathbb Q}^{o}(s):=\Big(\xi(2)\cdot\xi(s+2)\xi(2s+2)\Big)\cdot Z_{\mathbb Q}^{Sp(4)}(s).
\]
Consequently,
\[
\begin{aligned}
\xi_{Sp(4),\mathbb Q}^o(b)&\\
=&\frac{1}{b-1}\xi(2)\cdot\xi(b+2)\xi(2b+2)-\frac{1}{b+2}\xi(2)\cdot\xi(b)\xi(2b+1)\\
&-\frac{1}{2b}\cdot \xi(b+2)\xi(2b+2)+\frac{1}{2(b+1)}\cdot\xi(b) \xi(2b+1)\\
&-\frac{1}{(2b)(b+2)}\cdot\xi(b+1)\xi(2b+2)+\frac{1}{(-2b-2)(b-1)}\cdot\xi(2b+1)\xi(b+1),
\end{aligned}
\]
and that $\xi_{Sp(4),\mathbb Q}^o(b)$ satisfies the following functional equation
\[
\xi_{Sp(4),\mathbb Q}^o(-b-1)=\xi_{Sp(4),\mathbb Q}^o(b).
\]

\noindent
{\bf Definition \& Proposition.}
{\it Define the zeta function $\xi_{Sp(4),\mathbb Q}(s)$ for $Sp(4)$ over $\mathbb Q$ by
\[
\xi_{Sp(4),\mathbb Q}(s):=\xi_{Sp(4),\mathbb Q}^{o}(s-1).
\]
Then}

\noindent
(1) {\it The zeta function $\xi_{Sp(4),\mathbb Q}(s)$ is given by
\[
\begin{aligned}
\xi_{Sp(4),\mathbb Q}(s)
=&\frac{1}{s-2}\xi(2)\cdot\xi(s+1)\xi(2s)-\frac{1}{s+1}\xi(2)\cdot\xi(s-1)\xi(2s-1) \\
 &-\frac{1}{2s-2}\cdot \xi(s+1)\xi(2s)+\frac{1}{2s}\cdot\xi(s-1) \xi(2s-1) \\
 &-\frac{1}{(2s-2)(s+1)}\cdot\xi(s)\xi(2s)-\frac{1}{(2s)(s-2)}\cdot\xi(s)\xi(2s-1).
\end{aligned}
\]}

\noindent
(2) {\it It satisfies the standard functional equation
\[
\xi_{Sp(4),\mathbb Q}(1-s)=\xi_{Sp(4),\mathbb Q}(s).
\]}

\noindent
(3) {\it There are only four singularities, all simple poles, at $s=-1, 0, 1, 2$
and its residue at $s=2$ coincides with the volume of the compact domain $\frak F(0)$, that is,
\[
\mathrm{Res}_{s=2}\xi_{Sp(4),\mathbb Q}(s)=\xi(2)\xi(4)-\frac{1}{4}\xi(2)-\frac{1}{3}\xi(2)+\frac{1}{4}.
\]}

\noindent
{\it Remark.} The formula for the volume is arranged in the form to reflect
the fact that $\frak F(0)$ is obtained from the total fundamental domain $\frak F(\infty)$ whose
volume is given by $\xi(2)\xi(4)$, a result due to Siegel, by substracting two cuspidal
neighborhoods corresponding to two maximal parabolic subgroups whose volumes are
$\frac{1}{4}\xi(2)$ and $\frac{1}{3}\xi(2)$ respectively and adding a cuspidal neighborhood corresponding to
Borel subgroups whose volume is simple $\frac{1}{4}$. For related results, please refer to \cite{KW}.
\vskip 0.30cm
Furthermore, we have the following result of M. Suzuki (\cite{S2}):
\vskip 0.30cm
\noindent
{\bf Riemann Hypothesis for $\xi_{Sp(4),\mathbb Q}(s)$.}
\vskip 0.20cm
{\it All zeros of $\xi_{Sp(4),\mathbb Q}(s)$ lies on the central line $\mathrm{Re}(s)=\frac{1}{2}.$}

\subsubsection{Zetas for $Sp(2n)$}
Recall that  we have introduced the associated period $Z_{\mathbb Q}^{Sp(2n)}(z_n)$ for $Sp(2n)$ over $\mathbb Q$ by
\[
\begin{aligned}
Z_{\mathbb Q}^{Sp(2n)}(z_n):=
& \mathrm{Res}_{z_{n-1}-z_n=1}\cdots \mathrm{Res}_{z_{2}-z_3=1}\mathrm{Res}_{z_{1}-z_2=1} \\
& \Bigg(\sum_{w\in W} \frac {1}{\prod_{\alpha\in \Delta_0} \langle
  w\lambda-\rho, \alpha^{\vee}\rangle}  \cdot \prod_{\alpha>0, w\alpha<0}
  \frac {\xi\big(\langle\lambda, \alpha^{\vee}\rangle\big)}{\xi\big(\langle \lambda, \alpha^{\vee}\rangle+1\big)}\Bigg).
\end{aligned}
\]
Clearly,  there are minimal integers $I,\, J$, constants $a_i,\,b_i,\, i=1,2,\dots,I$, and $c_j, j=1,2,\dots,J,$
such that after multiplying the factor $\prod_{i=1}^I\xi(a_iz_n+b_i)\cdot\prod_{j=1}^J\xi(c_j)$,
the resulting function
\[
\xi_{Sp(2n),\mathbb Q}^o(z_n):=\Big(\prod_{i=1}^I\xi(a_iz_n+b_i)\cdot\prod_{j=1}^J\xi(c_j)\Big)\cdot Z_{\mathbb Q}^{Sp(2n)}(z_n)
\]
admits only finitely many poles and there are no special $\xi$-values appeared in the denominators.
\vskip 0.30cm
\noindent
{\bf Conjecture.} (Functional Equation) {\it There exists a constant $c_n$ such that
\[
\xi_{Sp(2n),\mathbb Q}^o(c_n-s)=\xi_{Sp(2n),\mathbb Q}^o(s).
\]}

\noindent
{\bf Definition.} {\it The zeta function $\xi_{Sp(2n),\mathbb Q}(s)$ of $Sp(2n)$ over $\mathbb Q$ is defined by
\[
\xi_{Sp(2n),\mathbb Q}(s):=\xi_{Sp(2n),\mathbb Q}^o\Big(s+\frac{c_n-1}{2}\Big).
\]}

\noindent
{\it Remark.} This in fact is a special case of a more general construction. For details, please see \cite{W4}.

As a direct consequence of the above conjecture, $\xi_{Sp(2n),\mathbb Q}(s)$
satisfies the following functional equation
\[
\xi_{Sp(2n),\mathbb Q}(1-s)=\xi_{Sp(2n),\mathbb Q}(s).
\]
With such a normalization, we then expect the following generalized RH for our zetas $\xi_{Sp(2n),\mathbb Q}(s)$.
\vskip 0.30cm
\noindent
{\bf The Riemann Hypothesis}$_{Sp(2n),\mathbb Q}$. \newline

{\it  All zeros of
$\xi_{Sp(2n),\mathbb Q}(s)$  lie on the line
$\mathrm{Re}(s)=\frac{1}{2}.$}

\noindent
Lin Weng \\
Graduate School of Mathematics \\
Kyushu University \\
6-10-1, Hakozaki, Higashi-ku \\
Fukuoka 812-8581, Japan \\
\texttt{weng@math.kyushu-u.ac.jp}
\smallskip

\noindent
{\it and}
\smallskip

\noindent
Chennai Mathematical Institute\\
Plot H1, SIPCOT IT Park\\
Padur PO, Siruseri 603103, India
\bigskip

\end{document}